\documentstyle[11pt]{article}
\newtheorem{thm}{Theorem}[section]
\newtheorem{prop}[thm]{Proposition}
\newtheorem{cor}[thm]{Corollary}
\newtheorem{lemma}[thm]{Lemma}
\newtheorem{rmk}[thm]{Remark}

\newcommand{\square}{\framebox(8,8){}}
\newenvironment{proof}{{\noindent\bf
Proof}\quad}{{\hfill\square}\vspace{4pt} \ \par}

\newenvironment{claim}{\vspace{5pt} \noindent{\bf
Claim}:\,}{\ \\ \par}

\newenvironment{example}{\addtocounter{thm}{1}
\vspace{5pt} \noindent{\bf Example \thethm}:\,}{\ \\ \par}

\renewcommand{\P}[1]{{\bf P}({#1})}             
\renewcommand{\sp}[1]{\mbox{${\bf P}^{#1}$}}    
\newcommand{\G}[1]{{\bf G}(1,{#1})}           
\newcommand{\C}{\mbox{$\bf C$}}     
\newcommand{\Q}{\mbox{$\bf Q$}}     

\newcommand{\pp}[1]{\mbox{$\scriptscriptstyle {#1}$}}

\title{On threefolds  covered by  lines }

\author{by Emilia {\sc Mezzetti} and Dario {\sc Portelli}\thanks{
This  work has been done in the framework of the
activities of EAGER. Both authors have been
supported by funds of MURST, Progetto Geometria Algebrica, Algebra
Commutativa e Aspetti Computazionali. The paper has been widely
improved while the second author was guest at the University of
Oslo. The financial support from both the University of Oslo and
the Norwegian Research Council is gratefully aknowledged.}}

\begin{document}

\date{}

\vskip 1 cm

\maketitle

\begin{abstract}
A classification theorem is given of projective threefolds 
that are covered by the lines of a two--dimensional family, 
but not by a higher dimensional family.
Precisely, if $X$ is such a threefold, let $\Sigma$ denote the Fano scheme
of lines on $X$ and  $\mu$  the
 number of lines contained in $X$ and passing through a general 
point of $X$. Assume that $\Sigma$ is generically reduced. Then $\mu\leq 6$. 
Moreover, $X$ is birationally a scroll over a surface ($\mu=1$),
or $X$ is a quadric bundle, or $X$ belongs to a finite list of 
threefolds of degree at most $6$. 
The smooth varieties  of the third type are precisely the Fano
threefolds with $-K_X=2H_X$.
\end{abstract}

\bigskip


\section*{Introduction}

\bigskip

\ \ \ \ Projective varieties containing \lq\lq\thinspace many'' 
linear spaces appear naturally in several occasions. For instance,
consider the following examples which, by the way, motivated our
interest in this topic.

The first example concerns  varieties of $4$-secant
lines of smooth threefolds in $\sp 5$. The family of such lines has in
general dimension four and the lines
fill up the whole ambient space, but it can happen that they form a
hypersurface.  

A second example comes from the following recent theorem of Arrondo 
(see \cite{A}), in some sense
the analogous of the Severi theorem about the Veronese surface:

{\it let $Y$ be a subvariety of dimension
$n$ of the Grassmannian $\G {2n+1}$ of lines of $\sp {2n+1}$ and assume
that $Y$ can be
isomorphically projected into $\G {n+1}$. Then, if the lines parametrized
by $Y$ fill up a variety
of dimension $n+1$, $Y$ is isomorphic to the second Veronese image of $\sp
n$.}

If those lines generate a variety of lower dimension, nothing is known.

\smallskip

In both cases it would be very interesting to have a classification 
of such varieties. Moreover, these examples show that for such a 
classification it would be desirable to avoid any assumption
concerning singularities. 

\smallskip

The first general  results about the classification of projective varieties
containing a higher
dimensional family of linear spaces were obtained by Beniamino Segre
(\cite{bS}). In particular,
in the case of lines, he proved:

{\it Let $X\subset\sp N$ be an irreducible variety of dimension $k$, let
$\Sigma\subset\G N$ be an
irreducible component of maximal dimension of the variety of lines
contained in $X$, such that the
lines of $\Sigma$ cover $X$. Then
$\dim\Sigma\leq 2k-2$. If equality holds, then $X=\sp k$. Moreover, if
$k\geq 2$ and
$\dim\Sigma = 2k-3$, then $X$ is either a quadric or a scroll in $\sp {k-1}$'s
over a curve.}

\smallskip

The case of a family $\Sigma$ of dimension $2k-4$ is treated in some
papers by Togliatti
(\cite{To}), Bompiani (\cite{Bo}), M. Baldassarri (\cite{B}), but their
arguments are not easy to
be followed. Recently, varieties of dimension $k\geq 3$
with a family of lines of dimension
$2k-4$ have been classified by Lanteri--Palleschi (\cite{LP}), as
particular case of a more general
classification theorem. Their starting point  is a pair $(X,L)$ where $L$
is an ample divisor on
$X$, which is assumed to be smooth or, more in general, normal and
$\Q$-Gorenstein. The
assumptions on the singularities of $X$ are removed  by Rogora in his
thesis (\cite {Ro}), but
he assumes $k\geq 4$ and codim $ X>2$.

\smallskip

The aim of this paper is the classification of the varieties of 
dimension $k$ covered by the lines of a family of dimension $2k-4$, 
in the first non--trivial case: $k=3$, i.e. threefolds 
covered by a family of lines of dimension $2$. So, we classify
threefolds covered by \lq\lq\thinspace few" lines.

\smallskip

A first remark is that among these varieties there are threefolds 
which are birationally scrolls over a surface or ruled by smooth 
quadrics over a curve. The first ones come from general surfaces
contained in $\G 4$, while the second ones come from general curves 
contained in the Hilbert scheme of quadric surfaces in $\sp n$. 
Note that these \lq\lq\thinspace quadric bundles" are built by 
varieties of lower dimension having a higher dimensional family of 
lines.

So we have focused our attention on threefolds not of these two types.

Observe that, if $X$ is a threefold covered by the lines of a family  
of dimension two, then there is a fixed finite number $\mu$
of lines passing through any general point of $X$. In particular, 
having excluded scrolls, we have assumed $\mu>1$.

\smallskip
It is interesting to remark that the surfaces $\Sigma$ in $\G 4$
corresponding to threefolds with $\mu>1$, can be characterized by the 
property that the tangent space to $\G 4$ at every point $r$ of 
$\Sigma$ intersects (improperly) $\Sigma$ along a curve. 
This follows from the fact that the points of $\G
4\cap T_{r}\G 4$ represent the lines meeting $r$.

Our point of view, that we have borrowed from the quoted paper of 
Mario Baldassarri, is the following. Since we do not care about
singularities, we are free to projected birationally into $\sp 4$ 
our threefolds to hypersurfaces of the same degree and with the same 
$\mu$. Hence, it is enough to classify  hypersurfaces in $\sp 4$ having
a family of lines with the requested properties.

If $X\subset\sp 4$ is a hypersurface of degree $n,$ then the equation
of $X$ is a global section $G\in\Gamma (\sp 4,{\cal O}_{\pp{\sp 4}}
(n)).$ The section $G$ induces in a canonical way a global section
$s\in\Gamma (\G 4,S^nQ),$ where $Q$ is the universal quotient bundle
on $\G 4.$ It
is a standard fact that the points of the scheme of the zeros of the 
section $s$ of $Q$ correspond exactly to the lines on $X.$ In this paper we will denote
by $\Sigma$ the  Fano scheme of the lines on $X$, which is, by definition, the scheme of the
zeros of the section
$s.$ 

In this paper we will study threefolds $X$ in $\sp 4$ covered by lines
such that $\Sigma$ has dimension two. 
\medskip

The following theorem is the main result of the paper.

\medskip

\begin{thm} \label{lista}
Let $X\subset\sp 4$ be a projective, integral hypersurface 
over an algebraically closed field $K,$ of characteristic zero,
covered by lines. Let $\Sigma$ denote the Fano scheme   of the lines on $X$ just
introduced. Assume that $\Sigma$ is generically reduced, 
 that $\mu>1$ and
that $X$ is not birationally ruled by quadrics over a curve. 
Then one of the following happens:
\begin{enumerate}

\item $X$ is a cubic hypersurface with singular  locus of dimension at most
one; if $X$ is smooth,
then $\Sigma$ is irreducible and $\mu=6$;

\item $X$ is a projection of a complete intersection of two hyperquadrics
in $\sp 5$; in general,
$\Sigma$ is irreducible and $\mu=4$;

\item $\deg X=5$: $X$ is a projection of a section of $\G 4$ with a $\sp
6$, $\Sigma$ is
irreducible and $\mu=3$;

\item $\deg X=6$: $X$ is a projection of a hyperplane section of $\sp
2\times\sp 2$, $\Sigma$ has
two irreducible components and $\mu=2$;

\item $\deg X\leq 6$: $X$ is a projection of  $\sp 1\times\sp 1\times\sp
1$, $\Sigma$ has
at least three irreducible components and $\mu\geq 3.$
\end{enumerate}
\end{thm}

Note that these five cases are precisely the projections of
Fano varieties with $-K_X=L\otimes L,$ \ $L$ ample
(\cite{LP}).
This list is the same as in the article of Baldassarri.

It is interesting to remark that the bound $\mu =6$ is attained  only
by cubic threefolds.

\medskip
The assumption that $\Sigma$ is generically reduced is necessary to make our method
work. Note that this is a genericity assumption for $X$ (however our threefolds are not
general, if the degree is $>3$; in fact, none of them is linearly normal in $\sp 4$, so
they have  a big singular locus).
This assumption  is quite strong, because it implies in
particular that the dual variety of $X$ is a hypersurface and that a general line on $X$ is
never contained in a fixed tangent plane.

The paper is organized as follows.

In \S \ \ref{par1} we prove that, under suitable conditions, 
on a general line of $\Sigma$ there are $n-3$ singular points 
of $X,$ where $n$ is the degree of $X,$ and we derive from this
many consequences we shall need in the paper. 
In particular, we will show 
that, if $n\geq 5,$ then the singular locus of $X$ is a 
surface and give an explicit lower bound for its degree
(Theorem \ref{delta}). 
Our main technical tool will be the family of
planes containing a line of $\Sigma$.
We prove that the assumption $\Sigma$ generically reduced implies that
there is no fixed tangent plane to $X$ along a general line on $X.$
From this it follows readily that the dual
variety of $X$ is a threefold (Theorem \ref{duale}). 
In this section we also introduce the ruled surfaces
$\sigma(r)$, generated
by the lines on $X$ meeting a fixed line $r$.

\S \ 2 contains the proof of the bound $\mu\leq 6.$ Moreover, if $n>3$
we prove that $\mu\leq 4.$

\S \ 3 is devoted to the classification of threefolds with an 
irreducible family of lines with $\mu>1$. First of all, we check
that, if $deg(X)>3$ and $X$ is not a quadric bundle,
only two possibilities are allowed for $\mu ,$ 
i.e. $\mu=3,4$. The threefolds with these invariants are then 
classified, respectively in Propositions \ref{mi3} and \ref{mi4}.

\S \ 4 contains the classification of threefolds with a reducible 
$2$-dimensional family of lines, such that all components of
$\Sigma$ have  $\mu_i=1$.

\bigskip

It is a pleasure to thank E. Arrondo, J.M. Landsberg, R. Piene
and K. Ranestad for several useful conversations about the
content of the paper, as well as for constant encouragement.

\bigskip

\noindent
In the paper we will use the following:

\medskip

\noindent {\bf Notations, general assumptions and conventions}

\begin{enumerate}

\item
We will always work over an algebraically closed field $K$ of
characteristic zero.

\item
$X\subset\sp 4$ will be a projective, integral hypersurface, 
of degree $n,$ covered by lines.

\item
We will denote by $\Sigma$  the Fano scheme 
of the lines on $X$.
(In particular, by the result of B. Segre quoted above, from 
$dim(\Sigma )=2$ it follows $n\geq 3.$)

\item
Let $\mu$ be the number of lines of $\Sigma$ passing
through a general point of $X$. If $\Sigma$  is reducible, with
$\Sigma_1,\ldots,\Sigma_s$ as irreducible components of dimension
$2,$ then we will denote by $\mu_i$ the number of lines of
$\Sigma_i$ passing through a general point of $X.$ Clearly
$\mu=\mu_1+\dots+\mu_s$. We assume $\mu >1.$

\item
We will assume that $X$ is not birationally ruled by quadric surfaces
over a curve.

\item
{}For a \lq\lq\thinspace general line in $\Sigma$'' we mean any line 
which belongs to a subset $S\subset\Sigma$ (never given explicitly),
such that $S$ is Zariski dense in $\Sigma.$ So \lq\lq\thinspace 
general line in $\Sigma$''  is meaningful also in the case of a 
reducible $\Sigma$.

\item
We will denote by the same letter both a line in $\sp 4$ and the 
corresponding point of $\G 4$. We hope that it will be always 
clear from the context which point of view is adopted.

\item
{}For $r$  general in $\Sigma$, the assumption $\mu>1$ ensures that 
the union of all the lines of $\Sigma$ meeting $r$ is a surface 
$\sigma(r)$, which can also be seen as a curve inside $\G 4.$
As $r$ varies in $\Sigma$, these curves  describe an
algebraic family in $\Sigma$ of dimension $\leq 2.$
If $\Sigma$  is reducible, with $\Sigma_1,\ldots,\Sigma_s$ as
irreducible components of dimension $2,$ then the
surfaces $\sigma(r)$ are unions $\sigma_1(r)\cup\ldots\cup
\sigma_s(r)$, where $\sigma_i(r)$ is formed by the lines of
$\Sigma_i$ intersecting $r.$ 

\end{enumerate}


\section{Preliminary results}\label{par1}


We consider  the degree $n$ of $X$.  For $n=3$,  it is
well known that all
cubic hypersurfaces of $\sp 4$ contain a family of
lines of dimension at least $2$, and of dimension exactly
$2$ if the singular locus of $X$ has codimension at least $2$.
For $n\geq 4$, a {\sl general} hypersurface of $\sp 4$ of degree $n$
is not covered by lines.

The following theorem is the main technical result of the paper. Here the assumption that
the irreducible components of dimension two of $\Sigma$ are reduced is essential.

\smallskip

\begin{thm}\label{n-3}
If $r$ is a general line of an irreducible component $\Sigma_1$
of $\Sigma$ which is of dimension two and generically reduced, 
then $r\cap Sing(X)$ is a
$0$-dimensional scheme of lenght $n-3$ (we will express this briefly
by saying that \lq\lq\thinspace on $r$ there are exactly
$n-3$ singular points of $X$"). 
In particular, if $n\geq 4$, then $X$ is singular.
\end{thm}

\begin{proof}
Let $r$ be a general line of $\Sigma_1$ (in particular, $\Sigma_1$ 
is the only component of $\Sigma$ containing $r$), and let $\pi$
be a plane containing $r$. Then $\pi\cap X$ splits as a union
$r\cup C$ where $C$ is a plane curve of degree $n-1$. So $r\cap C$ 
has length $n-1$ and it is formed by points that are 
singular for $\pi\cap X$, hence  either
tangency points of $\pi$ to $X$ or singular points of $X$.
We will prove that, if $\pi$ is general among the planes containing 
$r$, then exactly $n-3$ of these points are singular for $X$.
To this end, let us consider the family (possibily reducible) of 
planes ${\cal F}=\{\pi\vert \ \pi \supset r, r\in \Sigma \}$; 
its dimension is $4$.

\smallskip 

\noindent {\bf Claim}. {\sl The general plane through $r$ cannot
be tangent to $X$ in more that two points.}

\smallskip

\noindent {\sl Proof of the Claim}\ We have to prove that
$X$ does not possess a $4$-dimensional family 
of $k$-tangent planes, with $k>2.$ Assume by contradiction
that $X$ possesses such a  family ${\cal G}.$
Let $O$ be a general point of $\sp 4$, $O\not\in
X$. The projection $p_O:X\to \sp 3$, centered at $O$,  is a
covering of degree $n$, with branch locus a surface $\rho$ contained
in $\sp 3$. There is a $2$-dimensional
subfamily ${\cal G'}$ of ${\cal G}$ formed by the
planes passing through $O$: they project to
lines $k$-tangent the surface $\rho$.
Then $\rho$ satisfies the assumptions of the following lemma:

\begin{lemma}\label{tantetri}
Let $S\subset\sp 3$ be a reduced surface and assume that there
exists an irreducible subvariety $H\subset\G 3,$ with $dim(H)
\geq 2,$ whose general point represents a line in $\sp 3$ which
is tangent to $S$ at $k>2$ distinct points. Then $dim(H)=2$ and
$H$ is a plane parametrizing the lines contained in a fixed plane
$M\subset\sp 3,$ which is tangent to $S$ along a curve.
\end{lemma}

\noindent Therefore there
exists a plane $\tau$ tangent to $\rho$ along a curve of degree
$k$. But $\tau$ is the projection of a $3$-space $\alpha$
passing through $O$,
which must contain the planes of ${\cal G'}$. So these planes are
$k$-tangent also to $X\cap\alpha$, which
is a surface of $\sp 3$:  this means that
{\sl all} planes tangent to $X\cap\alpha$ are $k$-tangent. Since
$X\cap\alpha$ is not a plane, this is a contradiction.

\medskip 

\noindent Therefore, we have {\sl at least} $n-3$ singular points 
of $X$ on $r.$ Assume there are $n-2.$ 

\medskip

Let $H\subset\sp 4$  be a hyperplane containing $r.$ Let us 
denote by $\G H\simeq\G 3$ the Schubert cycle in $\G 4$ 
parametrizing lines contained in $H.$ Then, for 
general $H$ the intersection $\G H\cap\Sigma$ is proper, 
namely it is purely $0$-dimensional. In fact,
if infinitely many lines of $\Sigma$ were contained in
$H,$ then $dim(\Sigma )\geq 3,$ a contradiction.
Moreover, since we assume that $\Sigma$ is generically reduced,
both $\Sigma$ and $\G H$ are smooth at $r.$ 
We will show, now, that if $r$ contains $n-2$ singular points of $X,$
then $\Sigma$ and $\G H$ do not intersect
transversally at $r,$ and this will yield a contradiction. In fact,
$PGL(4)$ acts transitively on $\G 4,$ and we can use \cite{Kl} 
because we have assumed that our base field $K$ has characteristic 
zero.

\medskip

Before we start, let us recall briefly for the reader convenience
some basic facts about $T_r\G 4.$ 
Let $\Lambda \subset K^5$ be the $2$-dimensional linear subspace 
corresponding to $r,$ i.e. $r=\P\Lambda.$ 
Then $T_r\G 4$ can be identified with 
$Hom_K(\Lambda ,K^5/\Lambda ),$ hence for a non zero 
$\varphi\in T_r\G 4$ we have $rk\ \varphi =1$ or $2.$ 
In both cases we can associate to $\varphi$ in a canonical
way a double structure on $r.$ 
When $rk\ \varphi =1$ this structure is obtained
by doubling $r$ on the plane $\P{\Lambda\oplus Im(\varphi )},$ hence 
it has arithmetic genus zero (\cite{H2}). 
When $rk\ \varphi =2$ the doubling of $r$ is on a smooth quadric 
inside 
$\P{\Lambda\oplus Im(\varphi )}\simeq\sp 3,$ and the arithmetic
genus is $-1.$ In both cases we have $r\subset\P{\Lambda\oplus 
Im(\varphi)}$ and $\varphi\in T_r\G{\P{\Lambda\oplus Im(\varphi )}}.$

\medskip

To prove  the non transversality of $\Sigma$ and $\G H$ at $r,$
it is harmless to assume that $H$ is not tangent to $X$ at any
smooth point of $r.$ Therefore, the singularities of the surface
$S:=X\cap H$ on the line $r$ are exactly those points which are 
already singular for $X.$ 

To fix ideas, let $r$ be defined by the equations $x_2=x_3=x_4=0,$
$H$ defined by $x_4=0,$ and $S$ defined in $H$ by $\overline G=0.$
Then, the restriction to $r$ of the Gauss map of $S$ is given 
analytically as follows:

\vskip -0.2 cm

$$\alpha :P\mapsto [\overline G_{\pp{x_0}}(P),
\overline G_{\pp{x_1}}(P),\overline G_{\pp{x_2}}(P),
\overline G_{\pp{x_3}}(P)]. 
\label{gauss}
$$

\smallskip

\noindent
We can regard the $\overline G_{\pp{x_i}}(P)$'s as polynomials of 
degree $n-1$ in the coordinates of $P$ on $r.$ Since we assume $X$ has 
$n-2$ singular points on $r,$ the four polynomials 
$\overline G_{\pp{x_i}}(P)$ 
have a common factor of degree $n-2.$ Therefore, if we clean up 
this common factor, the above map can be represented analytically 
by polynomials of degree $1.$ Therefore, the double structure on 
$r$ defined by $\alpha$ has arithmetic genus $-1,$ and it 
arises from a {\sl non zero} vector $\varphi\in T_r\G H.$ 

Now, for every $P\in r$ which is a smooth point for $S$ we have 
$\alpha (P)= T_PS =T_PX\cap H,$ and in particular we have $\alpha (P)
\subset T_PX.$ This means that $\varphi$ is also a tangent vector to 
the Fano scheme $\Sigma$ of the lines on $X$ (see \cite{H}, pp. 
209-210), i.e. $\varphi\in T_r\Sigma $. Since we assume that $\Sigma_1$ is the only
component of
$\Sigma$ containing $r$ and $\Sigma_1$ is reduced at $r,$ by the usual criterion for
multiplicity one, we conclude that
$\G H$ and $\Sigma$ are not transversal at $r,$ and
the proof is complete (for general facts about intersections multiplicities
the reader is referred to \cite{Ful}).
\end{proof}

\begin{proof} {\sl of Lemma \ref{tantetri}}

\noindent  The lines in $\sp 3$ which
are tangent to $S$ are parametrized by a ruled threefold $K
\subset\G 3$: any line on $K$ corresponds to the pencil of
lines in $\sp 3$ which are tangent to $S$ at a fixed smooth point.
Then $H\subset K.$

$H$ is a surface: otherwise, a general point $O\in\sp 3$ would be 
contained
in infinitely many lines of $H$, therefore, every tangent line to a
general plane section $C$ of $S$ would be $k$-tangent to $C,$ with
$k>2,$ a contradiction.

Let $L\subset\sp 3$ be a line corresponding to a smooth point of $H;$
then $L$ is tangent to $S$ at least at points $P,Q,R$.
Since a general point of $K$
represents a line which is tangent to $S$ at a unique point, $K$ has
three branches at $L.$ We denote by $U_{\pp P}, U_{\pp Q}, U_{\pp R}$
the tangent spaces to these branches at 
$L,$ i.e. $U_{\pp P}\cup U_{\pp
Q}\cup U_{\pp R}$ is contained in the tangent cone to $K$ at $L.$ We 
have \ $U_{\pp P}\cap U_{\pp Q}=T_{\pp {L}}H.$
The intersection of this plane with $\G 3$ is the union of two
lines. Then, a direct, cumbersome computation proves
that these lines inside $\G 3$ represent respectively
{\sl the pencil of lines in $T_{\pp P}S$ through $Q$ and the
pencil of lines in $T_{\pp Q}S$ through $P.$}

\smallskip\noindent
{\bf Claim}. {\sl For a general point $L\in H$
we have $T_{\pp{L}}H\subset\G 3.$}

\smallskip\noindent
It is sufficient to show that $T_{\pp{L}}H\cap\G 3$ contains
three distinct lines.

{}From $U_{\pp P}\cap U_{\pp Q}\cap U_{\pp R}=T_{\pp {L}}H$ we 
get that the two lines of $T_{\pp{L}}H\cap\G 3$ are contained 
also in $U_{\pp R}.$
If we translate all this into equations, an easy computation shows 
that $T_{\pp P}S=T_{\pp R}S.$
By symmetry we get \ $T_{\pp P}S=T_{\pp Q}S=T_{\pp R}S.$ Therefore,
the three distinct lines in $\G 3$ which correspond to the pencils
in $T_{\pp P}S$ of centres respectively $P,Q,R$ are contained in
$T_{\pp{L}}H$, and the claim is proved.

\smallskip

By continuity, all the tangent planes  $T_{\pp{L}}H$ belong to
 one and the same
system of planes on $\G 3.$ Therefore, the tangent planes at two 
general
points of $H$ meet, and either $H$ is a Veronese
surface, or its linear span $\langle H\rangle$ is a $\sp 4.$

The first case is impossible
because the tangent planes to a Veronese surface fulfill a cubic
hypersurface in $\sp 5,$ whereas $\G 3$ is a quadric. On the 
other hand,
the quadric hypersurface $\G 3\cap\langle H\rangle$ in $\langle H
\rangle =\sp 4$ contains planes, hence it is singular. Therefore,
the hyperplane $\langle H\rangle$ is tangent to $\G 3$ at some point
$r$ and all the lines of $H$ meet the fixed line $r$
in $\sp 3.$ Were the lines of $H$ not lying on a unique plane through
$r,$ then {\sl any} plane $N$ through $r$ would contain 
infinitely many
lines $3$-tangent to the plane section $S\cap N$ of $S,$ a 
contradiction.
\end{proof}

\medskip

In the statement of Theorem \ref{n-3} we assume that an
irreducible component of $\Sigma$ is
generically reduced. We will give now a criterion that
leads to an easy way to check in practice if this hypothesis is
satisfied. 

We generalize a little and assume that an integral hypersurface 
$X\subset\sp N$ is covered by lines and that the dimension 
of the Fano scheme $\Sigma$ of lines on $X$ is $N-2.$ 
Let the line $r$ represent a general point of an irreducible 
component $\Sigma_1$ of $\Sigma ,$ of dimension $N-2,$ and
let $p$ be a general point of $r.$

Let $[x_0,\ldots ,x_N]$ be homogeneous coordinates in $\sp N.$
Assume that the line $r\subset X$  is
defined by $x_2=\ldots =x_N=0,$ and that the point $p$ is
$[1,0,\ldots ,0].$ We will work on the affine chart $p_{01}=1$
of the Grassmannian $\G N.$ Coordinates in this chart are
$p_{02},\ldots ,p_{0N},p_{12},\ldots ,p_{1N}$ and the line $r$
is represented by the origin. It is easy to see that a line $l$
in this affine chart contains the point $p$ if and only if its
coordinates satisfy the equations $p_{12}=\ldots =p_{1N}=0.$

Moreover, we will work on the affine chart $x_0=1$ of $\sp N,$
and we set $y_i:=x_i/x_0$ for $i=1,\ldots ,N.$ Then $p$ is the
origin. 

Let $G=G_1+G_2+\ldots +G_n=0$ be the 
equation of $X$ in this chart, where the $G_i$ are the 
homogeneous components of $G.$ We can assume that the tangent
space to $X$ at $p$ is defined by $y_N=0$, and we can consider
$y_1,\ldots ,y_{N-1}$ as homogeneous coordinates in ${\bf P}
(T_pX).$ Then, the line $r$ is represented in ${\bf P}(T_pX)$
by the point $[1,0,\ldots ,0].$ Finally, it is convenient to
write $G_i=F_i+y_NH_i,$ where the $F_i$'\thinspace s are 
polynomials in $y_1,\ldots ,y_{N-1}.$

\smallskip

\begin{prop}\label{criterio}
Assume that a hypersurface 
$X\subset\sp N$ is covered by lines and that the dimension 
of the Fano scheme $\Sigma$ of lines on $X$ is $N-2.$ 
Let the line $r$ represent a general point of an irreducible 
component $\Sigma_1$ of $\Sigma ,$ of dimension $N-2,$ and
let $p$ be a general point of $r.$ With the notations introduced 
above, $\Sigma_1$ is reduced at $r$ if and only 
if the intersection of the hypersurfaces in ${\bf P}(T_pX)$
defined by $F_i=0$ \ for $i=2,\ldots , n,$ is reduced 
at $[1,0,\ldots ,0].$ Or, equivalently, if the $(y_2,\ldots ,y_{N-1})$-primary
component of the ideal
$(F_2,\ldots ,F_n)\subset K[y_2,\ldots ,y_{n-1}]$ is $(y_2,\ldots ,y_{N-1})$. 
\end{prop} 

\begin{proof}
Let $s\subset\sp N$ be a line such that
$s\not\subset X$ \ and \ $p\in s.$ Let $A\subset\G N$ be
the Schubert variety parametrizing the lines in $\sp N$ which intersect
$s.$ The only singular point of $A$ is $s.$ In fact, it
is easily seen that $A$ is the intersection of $\G N$ with the
(projectivized) tangent space to $\G N$ at $s.$ In particular,
the points of $A$ different from $s$ are exactly the tangent
vectors to $\G N$ at $s$ which are of rank $1.$
Then, by using the facts on tangent vectors to Grassmannians
briefly recalled in the proof of Thm.\thinspace\ref{n-3},
it is easily seen that $A$ is the affine cone inside $T_s\G N,$ 
over a $\sp 1\times\sp {N-2}\subset{\bf P}(T_s\G N).$ 
It is clear that $\Sigma_1$ and $A$ intersect properly at $r.$

We claim that $\Sigma _1$ is reduced at $r$ if and only if $\Sigma _1\cap A$ is
reduced at $r$. Assume that
$\Sigma_1\cap A$ is reduced at
$r.$  Let $\cal O$ be the local ring of $\G N$ at $r$ and let $I$ and
$J$ denote respectively the ideals of $\Sigma_1$ and $A$ in 
$\cal O$. Then the Artinian ring ${\cal O}/I+J$ is reduced, i.e.
it is a field, and we want prove that
${\cal O}/I$ is reduced. The Cohen-Macaulay locus of $\Sigma_1$
is certainly open and non empty. So, by genericity, we can assume
that ${\cal O}/I$ is Cohen-Macaulay. We have $dim({\cal O}/I)=
N-2=ht(J).$ But $J$ is generated by a regular sequence of length
$N-2$ since $A$ is smooth at $r.$ Therefore, the same is true
for $J+I/I,$ being ${\cal O}/I$ a Cohen-Macaulay ring. But 
${\cal O}/I+J$ is a field, hence $J+I/I$ is the maximal ideal of 
${\cal O}/I$. It follows that this last ring is a regular local
ring.

Assume, conversely, that $\Sigma_1$ is  reduced at $r$.
Then $\Sigma_1\cap A$ is reduced at $r$ bevause $r$ is general
 and because of Kleiman'\thinspace s
criterion of transversality of the generic translate, already
used in the proof of Thm.\thinspace\ref{n-3}.

\smallskip

Denote by $B$ the Schubert 
cycle in $\G N$ parametrizing the lines in $\sp N$ through $p.$
A moment's thought shows that the local rings at $r$ of 
$\Sigma_1\cap A$ and $\Sigma_1\cap B$ are the same. Then we are
reduced to compute the ideal of $\Sigma_1\cap B$ inside
${\cal O}={\cal O}_{\pp{\G N,r}}.$ 

To do this, we replace the parametric 
representation of a general line $l$ containing $p,$ namely 
$y_1=t,$ and $y_i=p_{0i}t$ \ ($i\geq 2,$ where $t$ varies in the 
base field $K$) in all the equations $G_i=0,$ for $i\geq 1.$
{}From $G_1(t,p_{02}t,\ldots ,p_{0N}t)=0$ we get simply $p_{0N}=0.$
Then, since the $F_i$ are homogeneous polynomials, the other
generators for  the ideal of $\Sigma_1\cap A$ at $r$ are the
$F_i(1,p_{02},\ldots ,p_{0,N-1})$ \ $i=2,\ldots ,n.$  
An obvious change of variables completes the proof.
\end{proof}

\smallskip

\begin{example}\label{example1}
Let $X$ be the variety of the secant lines of a rational normal 
quartic curve $\Gamma\subset\sp 4.$ It is well known that the 
degree of $X$ is $3.$ On $X$ we have two families of lines of 
dimension two, each covering $X.$ We denote by $\Sigma_1$ 
the family of the secant lines of $\Gamma .$ By Terracini's Lemma, these lines are
also the fibres of the Gauss map. Hence  $dim(\check X)=2$ and for the family
$\Sigma_1$  we have $\mu_1=1$. 

Since $deg(X)=3,$ the intersection of $X$ with its tangent space along $r$
is a cubic  surface
which is singular along $r,$ hence ruled. These new lines 
form the second family $\Sigma_2$.

With a suitable choice of coordinates, a concrete case of
such an
$X$ is given by the equation:

$$
y_4+y_1y_4-y_2^2-y_3^2-y_1y_2^2-2y_2y_3y_4-y_4^3=0\ ,
$$

\smallskip

\noindent and the line $r$ defined by $y_2=y_3=y_4=0$ is one of
the secant lines of $\Gamma .$ 
Now $F_2=y_2^2+y_3^2$ and  $F_3=y_1y_2^2.$ Then, the
curves $F_2=0$ and $F_3=0$ do not intersect transversally at 
$[1,0,0]$, and $\Sigma_1$ is not reduced at $r.$ In fact, on 
any line of $\Sigma_1$ there are two points of $Sing(X).$
This shows that {\sl the hypothesis \lq\lq\thinspace $\Sigma_1$ is 
generically reduced" in Theorem \ref{n-3} is essential.} 

Note also that the curves $F_2=0$ and  $F_3=0$
intersect outside $[1,0,0]$ transversally at two points.
These points represent two lines on $X$ through $p$, which
belong to $\Sigma_2$. Therefore
  $\mu_2=2.$
\end{example}

\medskip

The following proposition deals with a delicate point, namely
the possibility for a general line $r$ of $\Sigma$ to be contained
in a plane which is tangent to $X$ at any point of $r.$ 
\smallskip

\begin{prop}\label{piantg}
Let $X\subset\sp 4$ be an irreducible hypersurface covered by
the lines of a family of dimension $2$ such that 
 $\Sigma$ is generically reduced. Let $r\in\Sigma$ be
general. Then
there is no plane containing $r$
which is tangent to $X$ at any general point of $r$.
\end{prop}

\begin{proof} 
Assume by contradiction that there exists a plane $M$ such that
 $M\subset T_qX$ for every $q\in r\cap X_{sm}.$ 
 We  perform some local computations
and we use the same notations as in  Proposition \ref{criterio}. So, let ${\bf A}
^4$ be an affine chart in $\sp 4,$ with coordinates $y_1,\ldots
,y_4.$ Assume that the origin is a general point $p$ of $X,$ and 
that $T_pX$ is defined by $y_4=0.$ Let $r$ and $M$ be defined
respectively by $y_2=y_3=y_4=0$ and $y_3=y_4=0.$
Let $G=G_1+G_2+\ldots +G_n=0$ be the equation of $X$ in this chart.
We write also $G_i=F_i+y_4H_i,$ where the $F_i$'\thinspace s are 
homogeneous polynomials in $y_1,y_2,y_3.$
Since the line $r$ is represented in ${\bf P}(T_pX)$
by the point $[1,0,0],$ we have

$$
F_i=y_1^{i-1}A_{i,1}(y_2,y_3)+y_1^{i-2}A_{i,2}(y_2,y_3)+\ldots 
+A_{i,i}(y_2,y_3)\ ,
$$

\smallskip

\noindent where the $A_{i,j}$ are homogeneous polynomials of degree
$j,$ or zero.
 
Now, if we move the origin of our system of coordinates to the point
$q\in r$
by a change of coordinates of type $Y_1=y_1-t$ and $Y_i=y_i$ for
$i=2,3,4$ and $t\in K$ \ (hence\ $q=(t,0,0,0)$),
then in the new system of coordinates $X$ is defined by the equation

$$
\widetilde G_t(Y_1,\ldots ,Y_4)=G(Y_1+t,Y_2,Y_3,Y_4)=Y_4+
\sum_{i=2}^{n}\{\thinspace F_i(Y_1+t,Y_2,Y_3)+
Y_4H_i(Y_1+t,Y_2,Y_3,Y_4)\thinspace \}
$$ 
$$
=(1+f(t))Y_4+\sum_{i=2}^{n}t^{i-1} A_{i,1}(Y_2,Y_3)+\ H.O.T. \ ,
$$

\smallskip

\noindent where $f(t)\in K.$ Now, since $M\subset T_qX$ for every 
$q\in r\cap X_{sm},$ the above equation shows that, necessarily the 
linear term of $\widetilde G_t$ belongs to the ideal $(Y_3,Y_4)$ for 
every $t\in K.$ Therefore, the linear forms $A_{i,1}(Y_2,Y_3)$ 
are in the ideal $(Y_3)$ for every $i\geq 2.$ But in this
case the curves in $\P{T_pX}$ defined by $F_i=0$ are either
singular at $[1,0,0]$, or with tangent line $y_3=0$ at $[1,0,0].$ 
This contradicts Prop.\thinspace\ref{criterio}, and the proof is
complete. 
\end{proof}

\medskip
From Proposition \ref{piantg} we will deduce the following very useful corollaries.

Let $\gamma:X\cdots\rightarrow \check{\bf P}^4$ be the Gauss map, 
which is defined on the smooth locus $X_{sm}$ of $X$. The closure 
of the image is $\check  X$, the dual variety of $X.$ If $\dim 
\check X<3$, then the fibres of $\gamma$ are linear subvarieties of 
$X,$ and the tangent space to $X$ is constant along each fibre.

\begin{cor}\label{duale}
Let $X\subset\sp 4$ be an irreducible hypersurface covered by
the lines of a family of dimension $2$ such that $\Sigma$ is generically reduced.
Then the dual variety $\check X$ of $X$ is a hypersurface of
$\check{\bf P}^4.$
\end{cor}
\begin{proof}
{}First of all, the dimension of $\check X$ must be at least $2:$
otherwise $X$ would contain a $1$-dimensional family of planes,
hence  a $3$-dimensional family of lines, a contradiction.
So assume by contradiction that $dim(\check X)=2.$ But then along each fibre of the Gauss
map there is  even a fixed tangent hyperplane, contradicting Proposition \ref{piantg}.
\end{proof}

\medskip
\begin{cor}\label{viapiano}
Let $X\subset\sp 4$ be an irreducible hypersurface covered by
the lines of a family of dimension $2$ such that $\Sigma$ is generically reduced.
Let $\Sigma_1$ be an irreducible component of $\Sigma$ of dimension two, such
that $\mu_1>1$. Let $r\in\Sigma_1$ be general, and set $\sigma_1(r)=\{
r'\in\Sigma_1\thinspace \vert\thinspace r\cap r'\neq\emptyset \}$. Then
$r\notin\sigma_1(r)$.
\end{cor}
\begin{proof}
Assume the contrary. Then, when $r'\in\sigma_1(r)$ moves on
$\sigma_1(r)$ to $r,$ the plane $\langle r'\cup r\rangle$
moves to a limit plane $M.$ The intersection $X\cap M$ is a
curve which has the line $r$ as a \lq\lq\thinspace double
component"\ ; in particular, this curve is singular along $r.$

Then $M\subset T_qX$ for every $q\in r\cap X_{sm}.$ In fact, 
if $M\not\subset T_qX,$ then $X\cap M$ would be smooth at $q,$
contradiction.
\end{proof}

\bigskip

Let $\cal F$ be the $4$-dimensional family of planes introduced 
in the proof of Theorem \ref{n-3}. We will 
consider now its subfamily $\cal F'$ of dimension $3$, formed by 
the planes generated by pairs of coplanar lines of $\Sigma$.

\begin{prop} \label{nota}
Let $\pi$ be a general plane of  $\cal F'$  generated by the 
lines $r$ and $r'$ of $\Sigma$. Then $\pi$ is tangent 
to $X$ at exactly $3$ points of $r\cup r'$ (but maybe $\pi$
is tangent to $X$ elsewhere, outside $r\cup r'$).
\end{prop}

\begin{proof} By Theorem \ref{n-3} there are two tangency points 
of $\pi$ to $X$ on $r$ and two on $r'$. The point $r\cap r'$ 
is singular for $X\cap \pi$, but it cannot be singular for $X,$ 
because, otherwise, 
letting $r$ and $r'$ vary, every point of $X$ would be singular.
So $r\cap r'$ is a tangency point of $\pi$ to $X$. Hence, $\pi$ is 
tangent to $X$ at exactly three points {\sl lying on $r$ or $r'$}.
\end{proof}

\medskip

To prove the next proposition, and also in the sequel,
we will need the following refined form of the connectedness
principle of Zariski, due to A.Nobile (\cite{N}):

\begin{lemma} \label{nobile}
Let $f:X\to T$ be a flat family of projective
curves, parametrized
by a quasi--projective smooth curve, such that the fibres $X_t$
are all reduced and $X_t$ is irreducible for  $t\neq 0$. 
Assume that, for $t\neq 0$,  $X_t$ has a fixed
number $d$ of singular points $P_1^{t},\ldots,P_d^t$ and 
that there exist $d$ sections $s_j: T\to X$ such that
$s_j(t)=P_j^t$ if $t\neq 0$, that $s_i(t)\neq s_j(t)$ if $i\neq j$ 
and
that $\delta(X_t, P_j^t)$ is constant (where $\delta(X_t, P_j^t)$
denotes the length of the quotient $\overline A/A$, $A$ being the
local ring of $X_t$ at $P_j^t$ and $\overline A$ its normalization).
 If the singularities of $X_0$ are
$s_1(0),\ldots,s_d(0),Q_1,\ldots,Q_r$,
then $X_0\setminus\{s_1(0),\ldots,s_d(0)\}$
is connected.
\end{lemma}

\begin{prop}\label{notre} Let $n\geq 4$ and let $\pi$ be 
a general plane of an arbitrary irreducible component of 
$\cal F'$. Then $\pi$ does not contain three lines of $\Sigma$.
\end{prop}

\begin{proof}
Assume by contradiction that $\pi$ contains the lines $r,r',r''$. 
Then the residual curve of $r$ in $\pi\cap X$ splits as $r'\cup 
r''\cup C$. Hence, by Lemma \ref{nobile}, there is a new tangency 
point on $r'\cup r''$, against Proposition \ref{nota}.
\end{proof}

\bigskip

Since our hypersurfaces $X\subset\sp 4$ contain \lq\lq\thinspace too
many" lines if $n\geq 4,$ it is quite natural that they are far from
general in the linear system of all hypersurfaces of $\sp 4$ of a fixed 
degree $n.$ In fact, it will turn out that, if $n\geq 4$ none
of them is linearly normal. Hence their singular loci have always
dimension $2.$ We will prove, now, directly this last property,
under the more restrictive assumption that $n\geq 5,$ which is
sufficient for our application of the theorem.

\smallskip

\begin{thm}\label{delta}
Let $X\subset\sp 4$ be a  hypersurface of degree $n\geq 5,$ covered 
by a family of lines $\Sigma$ of dimension $2$, with $\mu>1$ Let
$\Delta$ denote the singular locus of $X$. Then $\Delta$ is a surface.
If $X$ is not birationally ruled by quadrics, then
$\deg (\Delta) \geq 2(n-3)$.
\end{thm}

\begin{proof}
We assume by contradiction that $\Delta$ is a curve. Then every point 
of $\Delta$ belongs to infinitely many lines of $\Sigma$.
The curve $\Delta$ is not a line because every line  of 
$\Sigma$ meets $\Delta$ in $n-3$ points, and $n\geq 5.$
If $x\in X$ is general, from $\mu >1$ it follows that through $x$
there are two secant lines of $\Delta ,$ say $r$ and $s.$ 
By Terracini'\thinspace s lemma the tangent space to $X$ must 
be constant along $r$ and also along $s.$ Therefore, the plane
spanned by $r$ and $s$ is (contained in) a fibre of the Gauss map, 
hence it is contained in $X.$ So, through a general point of $X$ 
there is a plane on $X,$ contradiction.
This proves that $\Delta$ is a surface.

To prove the assertion on the degree, we consider a general
plane of $\cal F'$. If it intersects properly $\Delta ,$ then
this intersection contains at least $2(n-3)$ points, and the claim follows.
If the intersection is not proper, then $\Delta$ contains a family
of plane curves of dimension $3,$ hence it is a plane. Let $H$ be a
hyperplane containing $\Delta ;$ then $X\cap H$ splits as the union
of $\Delta$ with a surface $S .$ If $P\in S$ is general, there are at least
two lines on $X$ passing through $P.$ Each of them meets $\Delta ,$
hence is contained in $H,$ and therefore in $S.$ This shows that $S$
is a union of smooth quadrics.
\end{proof}

\medskip

We will give in the next proposition some generalities on the surfaces 
$\sigma(r)$.

\smallskip

\begin{prop}\label{sigma}
Let $X\subset\sp 4$ be a hypersurface of degree $n$ covered by
the lines of the family $\Sigma$ of dimension $2$, with $\mu\geq 2$.
Let $r$ be a general line of $\Sigma$ and $\sigma(r)$ be the
union of the lines of $\Sigma$ intersecting $r$. Then:

(i) $\sigma(r)$ is a ruled surface, having $r$ as line of 
multiplicity $\mu -1$;

(ii) if the surfaces $\sigma(r)$ describe, as $r$ varies in 
$\Sigma$, an algebraic family in $X$ of dimension $<2$, then $X$ 
is covered  by a $1$-dimensional family of quadrics such that 
there is one and only one quadric of the family passing through
any general point of $X$.
\end{prop}

\begin{proof}
The first assertion of $(i)$ is clear. To prove the second, it is 
enough to
observe that exactly $\mu -1$ lines of $\Sigma$, different from $r$,
pass through a general point of $r$, and that these lines are
separated by
the blow-up of $X$ along $r$.

The assumption of $(ii)$ means that, for every $r$, the lines of
$\Sigma$ intersecting $r$ intersect also infinitely many
other lines of the family, so $\sigma(r)$ is doubly ruled, hence it
is a smooth quadric, or a finite union of smooth quadrics. In the second
case, the algebraic family 
described by the surfaces $\sigma(r)$ has dimension two, so this
case is excluded.
\end{proof}

\medskip

We will refer to threefolds $X$ as in $(ii)$ as 
\lq\lq\thinspace quadric bundles''.

In the following we will analyze the self--intersection of the curves $\sigma(r)$ on
$\Sigma$ assuming it  positive. If the family of these curves is one--dimensional, then the
self--intersection is zero and $X$ is a quadric bundle. This is the reason why we exclude
quadric bundles in our classification.

Our final task concerning the surfaces
$\sigma (r)$ will be the determination of their degree. For this we need another
proposition.

\medskip

Let $r$ and $r'$ denote two general  lines in the same irreducible
component $\Sigma_i$ of $\Sigma .$ We will call $\overline\mu_i$ 
the number of lines of {\sl all} $\Sigma$ intersecting both
$r$ and $r'.$ 

Recall that, for every $r\in\Sigma ,$ the curve $\sigma(r)\subset
\Sigma$ (we switch our point of view, now) parametrizes the lines 
of $\Sigma$ intersecting $r.$ 
If $\Sigma$  is reducible, with $\Sigma_1,\ldots,\Sigma_s$ as
irreducible components of dimension $2,$ then the
curves $\sigma(r)$ are unions $\sigma_1(r)\cup\ldots\cup
\sigma_s(r)$, where $\sigma_i(r)$ is formed by the lines of
$\Sigma_i$ intersecting $r.$ Note that, if $\mu_i=1$ for some
index $i$ and $r\in\Sigma_i$, then $\sigma_i(r)$ is empty.

Then $\overline\mu_i$ is the intersection number $\sigma(r)\cdot
\sigma(r')$ on (a normalization of) $\Sigma.$

\smallskip

\begin{prop}\label{nocubic}
Let $X$ be a threefold such that $deg(X)>3.$ Let 
$r$ and $r'$ be two general lines in the same irreducible
component $\Sigma_i$ of $\Sigma .$
Then $\overline\mu_i =\mu -2$ (independent of $i$\thinspace !) 
\end{prop}

\begin{proof}
To evaluate $\overline\mu =\sigma(r)\cdot\sigma(r')$
we choose the lines $r$ and $r'$ so that they intersects at a point 
$p,$ smooth for $X.$ Since $deg(X)>3,$ by Proposition \ref{notre} we can 
also assume that $r$ and $r'$ are the only lines of $\Sigma$ contained
in the plane $\langle r\cup r'\rangle$, so that the lines 
intersecting both $r$ and $r'$ are those passing through $p$.
The conclusion follows from Corollary\thinspace\ref{viapiano}
\end{proof}

\medskip

\begin{prop} \label{grado}
Assume $\deg X\geq 4$ and let $r$ be a general 
line on $X.$ Then $\deg\sigma(r)=3\mu-4$.
\end{prop}

\begin{proof}
Note first that $\deg\sigma(r)$ is equal to the degree of the curve,
intersection of $\sigma(r)$ with a hyperplane $H$. We can assume 
$r\subset H;$ then $H\cap\sigma(r)$ splits in the union of 
$r$ with $m$ other  lines meeting $r$. Indeed, if $P\in H\cap
\sigma(r)$ and $P\not\in r$, there  exists a line passing through $P$ 
and meeting $r$, which is necessarily contained in $H$. 
Moreover,  $\sigma(r)$ and $H$ meet along $r$ with intersection
multiplicity $\mu -1$ (Proposition \ref{sigma}). 
Therefore $\deg\sigma(r)=\mu-1+m$.

To compute $m$, the number of lines meeting $r$ and contained in 
a $3$-space $H$, we can assume that $H$ is tangent to $X$ at a 
point $P$ of $r$. In this case $H$ contains the $\mu-1$ lines
through $P$ different from $r$. 
To control the other $m-(\mu-1)$ lines, we use the following 
degeneration argument. 

Since $H$ is tangent to $X$ at $p\in r,$ the intersection
multiplicity of $\Sigma$ and $\G H$ at $r$ is $2$ (this will be proved in \S 2, Lemma
\ref{moltepl}). According to the so called\lq\lq dynamical 
interpretation of the multiplicity of intersection",
in any hyperplane $H'$ \lq\lq close" to $H$ (if we are working 
over \C\ this means:  in a suitable neighbourhood of $H$ for the 
Euclidean topology of $\check{\sp 4}$) there are  two 
distinct lines $g,g'\in\Sigma$ which both
have $r$ as limit position when $H'$ specializes to $H.$ Note that
the lines $g$ and $g'$ are skew, because otherwise $g\in\sigma(g'),$ 
which becomes $r\in\sigma(r)$ when $H'$ specializes to $H,$ 
a contradiction with Prop.\thinspace\ref{viapiano}.

Therefore, we can choose a family of $3$-spaces $H_t$, 
parametrized by a smooth curve $T$, such that $H_0=H$ and, for
general $t$, $H_t$ is generated by two skew lines $r_t$ and $r'_t$, 
having both $r$ as limit position for $t=0$.
The lines in $H$ meeting $r$ come from lines in $H_t$ meeting 
either $r_t$ or $r'_t$. In other words, the intersections
$\sigma(r_t)\cap H_t$ and $\sigma(r'_t)\cap H_t$ both move to
$\sigma(r)\cap H$. Therefore to preserve the degree of these 
intersections, the remaining lines intersecting $r$ have to come from 
the $\overline\mu$ lines of $H_t$ meeting both $r_t$ and $r'_t$. 
Note that, if $l$ is one of these \lq\lq\thinspace remaining" lines,
then the multiplicity of $l$ in $\Sigma\cap\G H$ is $1.$ In fact,
otherwise, $H$ would be tangent to $X$ at some point of $l;$ but $H$
is already tangent to $X$ at $p,$ and $p\notin l.$ 
We can conclude by the previous proposition that 
$m=\overline\mu +\mu-1=2\mu-3$.
\end{proof}


\section{Bounds for $\mu$}


It is well known that for a surface covered by the lines of a 
$1$-dimensional family, there are at most two lines through any
general point.
The following theorem is the analogous for  threefolds.

\begin{thm}\label{sei}
Let $X\subset\sp 4$ be a $3$-fold covered by lines. Assume that the Fano scheme $\Sigma$ 
is generically reduced and of dimension $2$.  Then $\mu\leq 6$.
\end{thm}

\begin{proof}
It was already remarked in the Introduction that for the degree 
$n$ of $X$ we have $n\geq 3$. 
Let $p$ be a general point of $X$ and fix a system of affine  
coordinates $y_1,\ldots,y_4$ such that $p=(0,\ldots,0)$.  
Let $G=G_1+\ldots+G_n$ be the equation of $X.$ As usual, we assume
that $T_pX$ is defined by $y_4=0,$ and, moreover, we write
$G_i=F_i+y_4H_i\ ,$ for $i\geq 2.$

The polynomials $F_2,\ldots,F_n$ define (if not zero) 
curves in the plane ${\bf P}(T_pX).$
In particular, $F_2=0$ is a conic $C_2,$ whose points represent 
tangent lines to $X$ having at $p$ a contact of order $>2$, 
and $F_3=0$ is a cubic $C_3;$ the points of $C_2\cap C_3$ 
represent the tangent lines to $X$ having at $p$ a contact of 
order $>3$, and so on. Clearly  the points of ${\bf P}(T_{p}X)$ 
corresponding to lines contained in $X$ 
are exactly those of $C_2\cap C_3\cap\dots\cap C_n.$

We have $F_2\neq 0$ at any general point of $X$ because, otherwise 
$X$ would be a hyperplane of $\sp 4.$ On the other hand, since 
$deg(X)\geq 3,$ at any general point of $X$ we have also that 
$F_3$ is not a multiple of $F_2$ (\cite{GH1}, Lemma (B.16)). 
In particular, we have $F_3\neq 0,$ and $C_2$ is not contained in $C_3.$

Now,  $dim(\check X)=3,$ so $C_2$ is an irreducible conic (see \cite{bS2} 
or \cite{GKZ}), and we are done.
\end{proof}

\begin{rmk}
{\rm Actually, it is possible to give a proof of Theorem \ref{sei} which is independent of
Theorem \ref{duale}, hence of the assumption that $\Sigma$ is generically reduced.}
\end{rmk}

\begin{lemma}\label{moltepl}
{}For general $H\in\check{\sp 4}$ the intersection $\Sigma\cap\G H$
is proper. Moreover, if $r\in\Sigma$ is general and $r\subset H,$
then the intersection multiplicity of $\Sigma$ and $\G H$ at $r$
is always $\leq 2$ and it is $1$ if and only if $H$ is not tangent 
to $X$ at any point of
$r\cap X_{sm}.$ 
\end{lemma}

\begin{proof} 
The first part of the statement was already shown in the proof of
Thm.\thinspace\ref{n-3}. 

Moreover, in the same proof we saw that, if $H$ is not tangent to 
$X$ at some point of $r,$ then $T_r\Sigma$ and $T_r\G H$ are
transversal inside $T_r\G 4.$ In fact, the GCD of the polynomials
$\overline G_{\pp{x_1}}$ in $(\ref{gauss})$ has degree exactly
$n-3.$ Hence the double structure on $r$ they define has arithmetic
genus $-2$ and does not represent any vector in $T_r\G H.$ Therefore
$T_r\Sigma\cap T_r\G H=(0)$ and the intersection is transversal.

So we have proved that $i(r)=1$
if and only if $H$ is not tangent to $X$ at any point of $r.$
Hence, we assume now that $H$ is tangent to $X$ at some point of $r.$   
To show that $i(r)\leq 2$ we perform some local computations.
Let $[x_0,\ldots ,x_4]$ be a system of homogeneous coordinates
in $\sp 4$ such that the line $r$ is defined by the equations $x_2=x_3=x_4=0.$
Let $H=T_PX,$ where $P=[0,1,0,0,0]$ and $H$ is defined by $x_4=0.$ Let $[p_{01},\ldots
,p_{34}]$ be the related Pl\"ucker coordinates. So
$r$ has coordinates $[1,0,\ldots ,0],$ hence $p_{01}\neq 0.$ 
We will restrict, from now on, to work in the affine chart $U_{01}$ of $\G 4$ 
given by $p_{01}\neq 0;$ coordinates in this chart are $p_{02},p_{03}, 
p_{04}, p_{12}, p_{13}$ and $p_{14}$.
The equations of $\G H$ inside $U_{01}$ are $p_{04}=p_{14}=0$.
Then the general point of a line $r\in U_{01}\cap\G H$ is
$[s,1,p_{02}-sp_{12},p_{03}-sp_{13},0]$. 

In a suitable system of coordinates, the equation of $X$  is of the
 form:

$$
F=x_2\Psi x_0^2+x_3\Psi x_0x_1+x_4
\Psi x_1^2 + ax_2^2+bx_2x_3+cx_2x_4+\ldots +fx_4^2+ $$
\begin{equation} + \hbox{terms of degree $>2$ in
$x_2,x_3,x_4$} \label{equaprep}
\end{equation}
\noindent where $\Psi, a,
\ldots, f\in  K[x_0,x_1]$ are forms of degree $n-3$ and 
 $n-2$ respectively. Here we have used 
 the condition $r\subset X$. Moreover the homogeneous part of degree $1$ in
$x_2,x_3,x_4$ of $F$ can be normalized in this way because there is no fixed tangent
plane to $X$ along $r$ (\cite{CG}).

{}From $P\notin Sing(X)$ it follows that the coefficient of $x_1^{n-3}$ in the
polynomial $\Psi$ is not zero, and we can set $\Psi =
x_1^{n-3}+\rho_1x_0x_1^{n-4}
+\rho_2x_0^2x_1^{n-5}+\ldots $ . The point $P$ is $(0,0,0,0)$ in the affine
chart
$x_1\neq 0$, and if we dehomogeneize $F$ w.r.t. $x_1$ we get

\begin{equation}
{}^aF=({}^aF)_1+({}^aF)_2+\ldots =x_4+x_0x_3+\rho_1x_0x_4+
V(0,1,x_2,x_3,x_4)+\ldots
\label{fbello}
\end{equation}

\noindent where  $V\colon = ax_2^2+bx_2x_3+cx_2x_4+\ldots +fx_4^2$.

 The condition $r\subset
X$ implies that
$F(s,1,p_{02}-sp_{12},p_{03}-sp_{13},0)$ is identically zero as a polynomial in $s$. If
we set
$F(s,1,p_{02}-sp_{12},p_{03}-sp_{13},0)=\alpha +\beta s+\gamma s^2+\delta s^3+\ldots
$, then we
can compute $\alpha , \beta ,\gamma , \delta$ from  (\ref{fbello}), and we get:

$$\matrix{
\alpha =\overline ap_{02}^2+\overline bp_{02}p_{03}+\overline dp_{03}^2\hfill\cr
\cr
\beta = p_{03}+ \
\hbox{terms of degree $>1$ in $p_{12},p_{13},p_{02},p_{03}$}\hfill\cr
\cr
\gamma =-p_{13}+p_{02}+\rho_1p_{03}+\
\hbox{terms of degree $>1$ in $p_{12},p_{13},p_{02},p_{03}$}\hfill\cr
\cr
\delta = -p_{12}+\rho_1(-p_{13}+p_{03})+\rho_2p_{03}+\
\hbox{terms of degree $>1$ in $p_{12},p_{13},p_{02},p_{03}$}\hfill\cr}
$$

\smallskip

\noindent where $\overline a,\overline b,\overline d$ are the constant
terms of the polynomials
$a(x_0,1),b(x_0,1)$ and $d(x_0,1)$ respectively. Note that $\alpha,
\beta,\gamma, \delta$ are some of the equations of $\G H\cap\Sigma.$

By setting $\beta =\gamma =\delta =0$ we define inside the four dimensional
affine space $H$ a curve which is smooth
at $(0,0,0,0),$ the point in $\G H$ which represents $r.$ The direction of
the tangent line to
this curve at $r$ is given by the vector $(\rho_1 ,-1,1,0).$

Assume by contradiction that $i(r)>2.$
Then this vector annihilates $\alpha$, and $\overline a=0.$ It follows that
$x_0$ divides
$a(x_0,x_1)$, hence $x_2^2a(x_0,1)$ does not give any contribution to
$({}^aF)_2$. Therefore,
the reduction modulo $x_4$ (the equation of $T_PX$ in $\sp 4$)
 of the polynomial $({}^aF)_2$ is $x_3(x_0+\overline bx_2+
\overline dx_3).$ This is the equation of the conic $C_2$ embedded in 
$\P {T_PX}$. But, since the dual variety of $X$ has dimension $3$, by Theorem \ref{duale}
this conic should be smooth (\cite{bS2}, \cite{GKZ}).
\end{proof}

\smallskip

\begin{thm}\label{nomu5}
Let $X$ be a threefold such that $deg(X)>3.$ Then 
 $\mu \leq 4.$
\end{thm}

\begin{proof}
We analyze in detail the case $\mu=5$. A similar proof can be given if $\mu=6$.
For a different proof of this last case, see \cite{Si}.

Let us recall that $\overline\mu=3$, so given $r,r'\in\Sigma$ general and skew,
there are three lines $a,b,c\in\Sigma$ meeting  both $r$ and $r'$.

\smallskip
The lines $a,b,c$ are pairwise skew, otherwise $r,r'$ would fail to be skew.
Since $\overline\mu =3,$ there exists a third line in $\Sigma$, besides
$r$ and $r'$, meeting  both $a$ and $b.$ The same conclusion holds for
the pairs $(a,c)$ and $(b,c).$

\smallskip

\begin{claim}
If $r$ and $r'$ are general lines of $\Sigma$, then the three lines of
$\Sigma$ constructed above starting from the pairs $(a,b)$,   $(a,c)$
 and $(b,c)$ are distinct.
\end{claim}

Assume the contrary. Then there exists a
unique line $s\in\Sigma$,  different from both $r$ and $r'$, which meets
$a,b,c$.  Note that
 all the six lines $a$, $b$, $c$,
$r$, $r'$, $s$ are
contained in  the
linear span of $r$ and $r'$.

We consider now a family of pairs of lines $\{(r_t, r'_t)\}$
on $X$, parametrized by a smooth quasi--projective curve
$T$, such that $r_t$ and $r'_t$ are disjoint for a general $t\in T,$,
while for $t=0$ the lines $r_0$ and $r'_0$ meet at a  point $P$, general on
$X$.
Therefore for
$t$ general $\alpha_t:=\langle r_t, r'_t\rangle$ is a $\sp 3$: we get a
family of $3$-spaces
whose limit position $\alpha_0$ is the tangent space $T_{\pp P}X$.

We can assume that the plane of $r_0$ and $r'_0$ does not contain other lines
of $\Sigma$ (because $n>3$). For general $t$, we have three lines $a_t$, $b_t$,
$c_t$, meeting  $r_t$ and $r'_t$, and a third line $s_t$,  meeting
$a_t$, $b_t$ and $c_t$, which exists by assumption.
For $t=0$, the lines $a_0$, $b_0$, $c_0$ still meet $r_0$ and $r'_0$,
and $s_0$ meets $a_0$, $b_0$ and $c_0$. Hence $a_0$, $b_0$, $c_0$
pass through $P$. By Corollary \ref{viapiano}, $s_0$ cannot coincide
with $a_0$, $b_0$ or $c_0$, therefore by the assumption $\overline\mu=3$ and
$\mu=5$, either $s_0=r_0$ or $s_0=r'_0$.

Assume  $s_0=r_0$ .

 By Lemma \ref{moltepl}, the intersection multiplicity of
 $\G {\alpha_0}$ and $\Sigma$ is two at each of the five points corresponding
to the lines $r_0$, $r'_0$, $a_0$, $b_0$, $c_0$, therefore, by the
dynamical interpretation of the intersection multiplicity, there exist
four more lines in $\alpha_t$ moving to $r'_0$, $a_0$, $b_0$, $c_0$
 respectively. Let $u_t$ be a line of $\alpha_t$, having $r'_0$ as limit
position:  by Corollary \ref{viapiano} $r'_t\cap u_t=\emptyset$.

Let us assume that $u_t\cap(a_t\cup b_t\cup c_t \cup r_t\cup
s_t)=\emptyset$. In this case,
from $\overline\mu=3$, it follows that there exist six lines in $\alpha_t$,
three of them
meeting both
$s_t$
 and $u_t$, three meeting both $r_t$ and $u_t$.

The limit position of each of these six lines passes through $P$: but in
this way
we get too many lines passing through $P$ in $T_{\pp P}X$, contradicting
  the \lq\lq\thinspace multiplicity two '' statement  of Lemma \ref{moltepl}.

Therefore $u_t$ meets either $r_t$ (or, symmetrically, $s_t$) or $a_t$ (or,
symmetrically,
$b_t$ or $c_t$).

Case $(i)$: $u_t\cap r_t\neq\emptyset$.

\noindent In this case $u_t\cap s_t=\emptyset$, otherwise we would have
four lines meeting
both
$r_t$ and $s_t$. Also $u_t\cap a_t=\emptyset$ (and analogously
$u_t\cap b_t$ and $u_t\cap c_t$), otherwise the three lines $r_t$, $a_t$
and $u_t$ would be
coplanar. Therefore there exist three lines meeting $u_t$ and $s_t$, two
more lines meeting
$u_t$ and $a_t$, two meeting  $u_t$ and $b_t$, two meeting  $u_t$ and
$c_t$: summing up,
we get nine new lines.

We get again a contradiction with Lemma \ref{moltepl}, because we have found $16$ lines
tending to lines of $T_{\pp P}X$ passing through $P$. We conclude that $u_t\cap
r_t=\emptyset$.

Case $(ii)$: $u_t\cap a_t\neq\emptyset$.

\noindent So, being $\overline\mu=3$, $u_t\cap b_t=u_t\cap c_t=\emptyset$.
In this case, we
can construct four new lines, two meeting $s_t$ and $u_t$ and two meeting
$r_t$ and $u_t$.
Summing up we have $11$ lines moving
to lines of $T_{\pp P}X$ passing through $P$: this  contradiction proves
the Claim.

 Hence, given $r$ and $r'$, general lines on $X$, there exist lines $a$,
$b$ and $c$ meeting
both of them, and two by two distinct lines $s_1$, $s_2$, $s_3$ meeting $a$
and $b$, $a$
and $c$, $b$ and $c$ respectively. Moreover:  $s_i\cap s_j=\emptyset$ for
$i\neq j$;
$r\cap s_i=r'\cap s_i=\emptyset$,  $\forall i$.

Using the assumption $\overline\mu=3$, we get the existence of six more lines:
$l$ meeting $r$ and $s_1$, $l'$ meeting $r'$ and $s_1$; $m$ meeting $r$ and
$s_2$,
$m'$ meeting $r'$ and $s_2$; $n$ meeting $r$ and $s_3$, $n'$ meeting $r'$
and $s_3$.
Altogether there is a configuration of $14$ lines obtained from $r$ and $r'$.

The first observation is that the  $s_i$'s tend to lines through $P$, but $s_1$
tends neither to $a_0$ nor to $b_0$, because $s_1$ meets $a$ and $b$.
Therefore there
are three possibilities, that we examine separately:
\begin{description}
 \item [(i)] $s_1\rightarrow r_0$;
in this case the lines tending to $r_0$ are only $r$ and $s_1$. Now we
consider $s_2$: there
are two subcases:
   \begin{itemize}
     \item $s_2\rightarrow r'_0$: hence $s_3\rightarrow a_0$. Since $l'$
meets both
$r'$ and $s_1$, then it moves either to $b_0$ or to $c_0$; similarly $m'$,
which
meets both $r$ and $s_2$, moves either to $b_0$ or to $c_0$, and also $n$
does the same.
This contradicts Lemma \ref{moltepl} and Corollary \ref{viapiano}.
     \item $s_2\rightarrow b_0$: then we consider $l'$, which moves either
to $a_0$ or to
$c_0$. If $l'\rightarrow a_0$: then $s_3$, which meets $b$ and $c$, goes to
$r_0'$;
$m'$, which meets $r$ and $s_2$, goes to $c_0$; $n$ which meets $r$ and $s_3$
could go to $a_0$ or to $b_0$ or to $c_0$: but all three cases are excluded by
Lemma \ref{moltepl} and Corollary \ref{viapiano} again. If  $l'\rightarrow c_0$, the
conclusion is similar.
    \end{itemize}
 \item [(ii)]$s_1\rightarrow r'_0$; this case is analogous to case (i).
 \item [(iii)]$s_1\rightarrow c_0$. We consider $s_2$: since it meets $a$
and $c$, it
goes to $b_0$, or to $r_0$, or to $r_0'$. The last two possibilities are
excluded as in (i)
and (ii) for $s_1$, so $s_2\rightarrow b_0$ and finally $s_3\rightarrow
a_0$. By considering
the limit positions of $l$, $l'$, $m$, we find that also in this case the
\lq\lq\thinspace
multiplicity two'' statement of Lemma \ref{moltepl} is violated.
\end{description}
This concludes the proof.
\end{proof}

\medskip

\smallskip

The statement of Theorem \ref{lista}
shows that the  families of lines in $\sp 4$ we want to
classify are characterized by  the number $s$ of irreducible components
$\Sigma_1$,
$\ldots$, $\Sigma_s$  of $\Sigma$ and by the relative $\mu_i$'s. Therefore the
proof can be  organized according to the following two possibilities:

\begin{itemize}

\item there exists an irreducible component $\Sigma_i$ of $\Sigma$ 
with $\mu_i>1$;

\item for every irreducible component $\Sigma_i$ of 
$\Sigma$,\  $\mu_i=1$.

\end{itemize}

\noindent
By Theorem \ref{sei}, there are only finitely many values of $s$ and $\mu_i$
to analyze. A posteriori, it will turn out that, actually, in the first 
case there do not exist other irreducible components of $\Sigma .$


\section{There exists an irreducible component $\Sigma_i$ of $\Sigma$ 
with $\mu_i>1$}


Let $\Sigma_i$ be an irreducible component of $\Sigma$ of dimension $2,$
such that $\mu_i>1.$ In this section we will consider and use only the lines
of $\Sigma_i$, e.g. for constructing the surfaces $\sigma (r)$ and so on.
So, for simplicity, we will denote $\Sigma_i$ by $\Sigma$ and $\mu_i$ by
$\mu.$ Note that Proposition \ref{nocubic} is still true (with the
same proof) even if we use in the statement our \lq\lq\thinspace$\mu$"
and \lq\lq\thinspace$\overline\mu$" defined by using only the lines
of $\Sigma_i.$

\smallskip

\begin{prop}\label{nomu2}
Assume that $X$ is not a quadric bundle and that $deg(X)>3.$
Then $\mu >2.$
\end{prop}

\begin{proof}
Since we assume that $X$ is not a quadric bundle we have that the
dimension of $\{\thinspace \sigma(r)\thinspace\}_{\pp{r\in\Sigma}}$
is $2$ by Prop.\thinspace\ref{sigma}. Then, through a general point
of $\Sigma$ there are infinitely many curves $\sigma (r),$ and,
by  Proposition \ref{nocubic} we conclude

$$
\mu -2 =\sigma (r)^2>0\ .
$$
\end{proof}

\medskip

 Then, if we assume that
$deg(X)>3$ and that  $X$ is not a quadric bundle, by the above proposition 
and by Theorem\thinspace\ref{nomu5}, the only possibilities for 
$\mu$ are $\mu =3,4.$

\bigskip

\noindent
{\bf The case $\mu =3.$}

\medskip

\begin{prop}\label{mi3}
Let $X\subset\sp 4$ be a hypersurface of degree $>3,$
containing an irreducible family of lines $\Sigma$ with $\mu =3.$
Then $X$ has degree $5,$ sectional genus $\pi=1$ and it 
is a projection of a Fano threefold of $\sp 6$ of the form $\G
4\cap\sp 6$.
\end{prop}

\begin{proof}
The algebraic system of dimension two
$\{\sigma(g)\}_{g\in\Sigma}$ on the surface $\Sigma$ is  linear because
there is exactly one curve of the system passing through two general points
($\overline\mu=1$).
Also the self--intersection is equal to $\overline\mu=1$,
therefore $\{\sigma(g)\}$ is a homaloidal net of rational curves, which
defines
 a birational  map $f$ from $\Sigma$ to the plane, such that the curves of
the net correspond
to the lines of $\sp 2$. The degree of the curves
$\sigma(g)$ is $5$ by Proposition \ref{grado}. So the birational inverse
of $f$
is given by a linear system of plane curves of
degree $5$. Hence we get immediately
the weak bound $\deg \Sigma\leq 25$.
Let $\nu$ denote the number of lines of $\Sigma$ contained in
a $3$-plane: by Schubert calculus, $\deg\Sigma=\mu n+\nu$.
To evaluate $\nu$, we consider two general skew lines $r$, $r'$ on $X$,
generating a $3$-space $H$. The lines $r$ and $r'$
have  a common secant line $l$. The set--theoretical intersection
$\sigma(r)\cap H$ is the union of
$r$, $l$ and two more lines $l_1$, $l_2$ by Proposition \ref{grado}.
Similarly we get two new lines $m_1$, $m_2$ in
$\sigma(r')\cap H$. The line $l_1$ (resp. $l_2$) cannot meet both $m_1$ and
$m_2$ because
$\overline\mu=1$, so there are two new lines in $H.$

So we have found at least $9$ lines in $H$, hence $\nu\geq 9$.
The assumption $\mu=3$ together with $\nu\geq 9$ gives at once
$n\leq 5$.

Let $S$ be a general hyperplane section of $X$. If $n=4$, then it 
is  well known (see for example the classical book of Conforto 
\cite{C}) that under our assumptions one of the following happens:
$S$ is a ruled surface (in particular a cone) or  a  Steiner surface 
or  a Del Pezzo surface with a double irreducible conic. None of these 
surfaces is section of a threefold $X$ with the required properties.
In the first case $X$ would have a family of lines of dimension
$3$, in the second case $X$ would be a cone, in the third case 
$\mu=4$ (see \cite{GH2} and \cite{T}). Therefore the degree of
$X$ is exactly $5.$

We can apply, now, Theorem \ref{delta} which gives   $\deg\Delta\geq
4$ since $n=5.$ If $\pi$ denotes the sectional genus of $X$ (i.e.
the geometric genus of a general plane section of $X$) we deduce 
$\pi\leq 2.$ 

To exclude $\pi=2,$ we show that there
exist planes containing
three lines of $\Sigma$. Indeed
let $r$ be a general line of $\Sigma$. We fix in
$\sp 4$ a $3$-plane $H$  not containing $r$,
intersecting $r$ at a point $O$. Let
$\gamma:=\sigma(r)\cap H$ be a hyperplane section of $\sigma(r)$.
By Proposition \ref{grado}, $\sigma(r)$ has degree $5$, hence
there exists a trisecant line $t$ passing through $O$ and meeting
$\gamma$ again at two points $P$ and $Q$.
Let $M$ be the plane generated
by $r$ and $t$: it contains also the lines of $\sigma(r)$ passing through
$P$ and $Q$, so $M$ contains three lines contained in $X$. Now we consider
$M\cap\Delta$.
By Lemma \ref{nobile} in $M\cap X=r\cup r'\cup r''\cup C$ there must be a
\lq\lq\thinspace new''
tangency point, hence  $\Delta\cap M$ contains at least five points.
Therefore $\deg\Delta\geq 5$
and $\pi\leq 1$.
If $\pi=0$,  the curves  intersection of
$S$ with its tangent planes have a new singular point, so they split. 
Then by the Kronecker--Castelnuovo 
theorem,  $S$ is ruled, a contradiction. So we have  $\pi=1$ and
$S$ is a projection of a linearly normal Del Pezzo surface $S'$ of $\sp 5$
of the same degree $5$
(see
\cite{C}), which is necessarily a linear section of $\G 4$.
This proves the theorem.
\end{proof}

\bigskip

\noindent
{\bf The case $\mu =4.$}

\medskip

\begin{prop}\label{mi4}
Let $X\subset\sp 4$ be a hypersurface of degree $>3,$
containing an irreducible family of lines $\Sigma$ with $\mu =4.$
Then $X$ has degree $4$ and sectional genus $\pi=1$,
hence it is a projection of a Del Pezzo threefold of $\sp 5,$  complete
intersection of two quadric hypersurfaces of $\sp 5$.
\end{prop}

\begin{proof}
Let ${\overline g}\in\Sigma$ be general and set $\sigma
\colon =\sigma({\overline g}),$ for simplicity. Let $\gamma$
denote a normalization of $\sigma.$
The proof of the proposition is based on the following two lemmas.

\begin{lemma}\label{curvabuona}
The curve $\gamma$ is irreducible, hyperelliptic of genus $2.$ 
Hence $\gamma$ can be embedded into $\sp 3$ as a smooth quintic.
\end{lemma}

Let $S\subset {\bf G}(1,3)$ be the surface parametrizing the secant
lines of $\gamma.$
Let $r\subset\sp 3$ be a fixed general secant line of $\gamma;$ we will
denote by $A$ and $B$ the points of $r\cap\gamma.$
The family of all secant lines of $\gamma$ that intersect $r$ has
three irreducible components: the secant lines through $A,$ those
through $B$ and \lq\lq the other ones". This last component is represented
on $S$ by an irreducible curve that we will denote by $I_r.$

\begin{lemma}\label{famgiusta}
There exists a birational map $\tau\colon\Sigma\cdots\to S$ such that
the image via $\tau$ of every curve $\sigma(g)\subset\Sigma$ is the
curve $I_{\tau (g)}$ on $S$ just introduced.
If $g,g'\in\Sigma$ are general, then $g\cap g'\neq\emptyset$
if and only if $\tau (g)\cap\tau (g')\neq\emptyset.$
\end{lemma}

We will prove now Proposition \ref{mi4} assuming Lemmas
\ref{curvabuona} and \ref{famgiusta} .

Let $p$ be a general point of $\sp 3,$ \  $p\notin\gamma.$ There
are four secant lines $l_1,\ldots ,l_4$ of $\gamma$ through $p$
and we can assume that $l_i=\tau (g_i),$ with $g_i\in
\Sigma ,$ \  $i=1,...,4.$
By Lemma \ref{famgiusta} we have $g_i\cap g_j\neq\emptyset$ for every
 $i\neq j.$

The first possibility is that, for a general $p\in\sp 3,$ the four
lines $g_1,\ldots ,g_4$ all lie in a plane $M_p\subset\sp 4.$
By Prop.\thinspace\ref{notre} the family of such planes has dimension
at most $2$ and,
therefore, the same plane $M_p$ corresponds to infinitely many points
of $\sp 3.$ This implies that every plane $M_p$ contains infinitely many
lines of $\Sigma,$ hence $M_p\subset X.$ Then $X$ contains at least
a $1$-dimensional family of planes: a contradiction.

Therefore, for a general $p\in\sp 3,$ the four lines $g_1,\ldots ,g_4$
all contain one fixed point $P\in X,$ and we get a rational
map $\alpha\colon\sp 3\cdots\to X$ by setting $\alpha (p):=P.$ This map
is dominant because $\tau\colon\Sigma\cdots\to S$ is birational, and it has
degree $1,$ because $\mu=4.$ Hence $X$ is birational to $\sp 3$ via $\alpha.$

\smallskip

Note that $\alpha$ is not regular at the points of $\gamma,$ so
$\alpha$ is defined by a linear system of surfaces $F\subset
\sp 3$ of degree $m,$  all containing $\gamma.$ Let $s$ be the
maximum integer such that these surfaces contain the $s^{th}$ infinitesimal
neighbourhood of $\gamma .$ So $F\in\vert mH-(s+1)\gamma\vert ,$ where
$H$ is a plane divisor in $\sp 3.$ We claim that $s=0$ and $m=3.$

The second part of the statement of Lemma \ref{famgiusta} makes clear
that {\sl any secant line of $\gamma$ is transformed by $\alpha$ into
a line of $\Sigma .$}
Therefore we must have $m=2(s+1)+1;$ if we intersect one of the surfaces $F$
with the unique quadric surface $Q$ containing $\gamma ,$ by Bezout
and $deg(\gamma )=5$ we get

$$
2m=2\big[ 2(s+1)+1\big] \geq 5(s+1),
$$

\noindent
hence $s\leq 1.$

If $s=1$ we get $m=5$ and the surfaces $F$ contain the first infinitesimal
neighbourhood of $\gamma .$
Let $I\subset K[x_0,\ldots ,x_3]$ denote  the saturated  ideal of $\gamma .$
Since $\gamma\subset\sp 3$ is arithmetically  Cohen-Macaulay, the saturated
ideal of the  first infinitesimal  neighbourhood of $\gamma$ is $I^2$
(\cite{Pe},
2.3.7). Now, $I$ can be minimally generated by one polynomial $q$ of
degree $2$ (the equation of $Q$) and two polynomials of degree $3;$
therefore, {\sl every homogeneous polynomial of degree $5$ in $I^2$
must contain $q$ as a factor.} So the case $s=1$ is excluded.

Hence, the linear system defining $\alpha$ is a system of {\sl cubic}
surfaces of $\sp 3,$ containing $\gamma$ with multiplicity $1.$ The linear
system  of all such surfaces  defines a rational
map $\sp 3\cdots\to\sp 5,$ whose image is a Del Pezzo threefold, complete
intersection of two quadric hypersurfaces of $\sp 5$. This completes the
proof of Proposition \ref{mi4}.
\end{proof}

\medskip

\begin{proof} {\bf of Lemma \ref{curvabuona}}
\ The proof is divided into several steps.

\smallskip
\noindent{\bf Step 1.}
\smallskip
\noindent
{\sl There is a birational map
$\psi\colon\Sigma\cdots\to\sigma^{(2)},$
where $\sigma^{(2)}$ denotes the symmetric product of the curve
$\sigma$ by itself.}

\smallskip

On  $\Sigma$ there is the algebraic system of curves
$\{\thinspace\sigma(g)\thinspace\}_{\pp{g\in\Sigma}}\ ,$  of
dimension $2.$
Since $\overline\mu =2$, there are exactly $2$
curves of the system containing two fixed general points on $\Sigma;$
moreover $\sigma(g)^2=2$.

\smallskip

The map $\psi$ is defined as follows:
let $r$ be a general line of $\Sigma$; let $a,b$ be the two lines of $\Sigma$
intersecting both $r$ and
${\overline g}.$ The corresponding points on $\Sigma$ actually lie on $\sigma.$
We set \ $\psi\colon r\mapsto (a,b);$ \ it is easily seen that $\psi$
is birational. Note that the map $\psi$ depends on the choice of
${\overline g}\in\Sigma.$

In particular, from $\Sigma$ irreducible it follows that $\sigma$
{\sl is also irreducible.}

\smallskip
\noindent{\bf Step 2.}
\smallskip
\noindent {\sl The characteristic series of the algebraic system
$\{\thinspace\sigma(g)\thinspace\}_{\pp{g}}$ on the curve $\sigma$
is a complete $g^1_2.$ Therefore also the algebraic system
$\{\thinspace\sigma(g)\thinspace\}_{\pp{g}}$ is complete.}

\smallskip

{}From the fact that the dimension and the degree of the algebraic system
$\{\thinspace\sigma(g)\thinspace\}_{\pp{g}}$ are both $2$, it
follows at once that the characteristic series has degree $2$ and
dimension $1,$ i.e. it is a $g^1_2.$

Assume it is not complete; then $\sigma$ is necessarily a rational
curve and the characteristic series generates a complete $g^2_2$.
 In this case $\Sigma$ is
a rational surface and we can embed $\{\thinspace\sigma(g)\thinspace\}_
{\pp{g}}$ into the complete linear system $\vert\sigma(g)\vert$ of
dimension $3$.
Let $L$ be the linear span of $\{\thinspace\sigma(g)\thinspace\}
_{\pp{g}}$ inside $\vert\sigma(g)\vert.$
Let $\cal L$ be the linear system of those ruled surfaces on
$X$ which correspond to the curves of $L.$
{}Fix a general point $P$ of $X$ and denote by $\cal M$ the subsystem of
surfaces of $\cal L$ containing $P$: $\cal M$ contains $4$ 
linearly independent
surfaces, hence its dimension is at least $3$: a contradiction.

\smallskip
\noindent{\bf Step 3.}
\smallskip
\noindent {\sl Let $\pi$ denote the geometric genus of $\sigma.$ Then
$\pi\geq 2.$}

\smallskip

By the previous step we already know that $\pi\geq 1;$ assume $\pi =1.$
Then, by the well known fact that the irregularity of $\sigma^{(2)}$
equals the (geometric) genus of $\sigma,$ the irregularity of $\Sigma$
is $1.$
But the surface $\Sigma$, which  parametrizes the curves of
$\{\thinspace\sigma(g)\thinspace\}_{\pp{g\in\Sigma}},$ is therefore
 fibered by a $1$-dimensional family of lines, each line representing a
linear pencil of curves $\sigma(g);$ from $\sigma(g)^2=2$ it
follows that every such pencil has $2$ base points.
This also means that on $X$ we have a $1$-dimensional family of linear
pencils of elliptic ruled surfaces $\sigma (g),$ each pencil having exactly
two base lines.

We fix one of these pencils $\{\thinspace\sigma (g_{\pp t})\thinspace\}_{t\in\sp
1},$ and we let $r$ and $r'$ denote the two base lines.
Every surface of the pencil is of the type $\sigma (g),$ with $g$ intersecting
both $r$ and $r'.$ Set

$$
R\colon = \bigcup_{t\in\sp 1}\thinspace g_t\subset X
$$

\noindent
We claim that, for general $t,t'\in\sp 1,$ the lines $g_t$ and $g_{t'}$
don't meet on $r.$ Indeed, if  $g_t\cap g_{t'}=P\in r$, then
also the fourth line of $\Sigma$ through $P$ would be contained in
$\sigma (g_{\pp t})\cap\sigma (g_{\pp t'}),$ the base locus of the pencil:
a contradiction.

So $r$ is a simple unisecant for $R$. Since $\sigma (r)$ is irreducible,
from $R\subseteq\sigma (r)$ it follows that $R=\sigma (r).$ Then we have
a contradiction because $r$ has multiplicity $3$ on $\sigma (r)$
by Proposition \ref{sigma}.
Therefore,  $\sigma$ is hyperelliptic of geometric genus
$\pi\geq 2$.

To complete the proof of Lemma \ref{curvabuona} it remains to show:

\smallskip
\noindent{\bf Step 4.}
\smallskip
\noindent
{\sl The genus of $\gamma$ is $2.$ In particular, $\gamma$ is embedded in
$\sp 3$ with degree $5.$}

\smallskip

Let $p\in\overline g$ be a general point, and let $a,b,c\in\Sigma$
denote the lines through $p,$ different from $\overline g.$ Moreover,
let $d,e\in\sigma$ be such that $d+e\in g^1_2$ on $\gamma .$ Then
$H:=a+b+c+d+e$ is a positive divisor on $\gamma ,$ of degree $5.$
When $p$ varies on $\overline g,$ the divisors on $\gamma$ of type
$a+b+c$ are all linearly equivalent because they are parametrized by
the rational variety $\overline g.$ We denote by $\cal D$ the pencil
of such divisors. Since the two rational maps $\gamma\to\sp 1$ defined
respectively
by $\cal D$ and $g^1_2$ are clearly different, it is easily seen that
$dim\ \vert H\vert\geq 3.$ Hence, by Clifford'\thinspace s theorem 
$H$ is non special. Since $\pi\geq 2,$ it follows then by Riemann--Roch
that $dim\ \vert H\vert = 3,$ and that $\pi =2.$ Then $H$ is also very
ample on $\gamma .$
\end{proof}

\medskip

To prove Lemma \ref{famgiusta} we need

\begin{lemma}\label{gradoind2}
 $\{\thinspace I_r\thinspace\}_{\pp{r\in S}}$\ is an algebraic
system of curves on $S$ of dimension $2,$ degree $2$ and index $2.$
\end{lemma}
\begin{proof}
Since $deg(\gamma )=5$ and $\pi =2,$ there are $4$ secant lines
of $\gamma$ through a general point of $\sp 3,$ and $10$ secant lines
of $\gamma$ contained in a general plane of $\sp 3.$ Therefore, the
class of $S$ in the Chow group $CH_2({\bf G}(1,3))$ is $4\alpha +10\beta ,$
with traditional notations. It follows that the degree of $S\subset\sp 5$
is $14;$ this means that there are $14$ secant lines of $\gamma$
intersecting two general lines $r$ and $r'$ in $\sp 3.$

Assume, now, that $r$ and $r'$ are chords of $\gamma ,$ and
set $r\cap\gamma =\{ A,B\} ,$ \ $r'\cap\gamma =\{ C,D\} .$
To compute $I_r\cdot I_{r'}$ we have just to compute the number of the
spurious solutions among these $14$ secant lines.
Let $M$ be the plane generated by $r$ and $C;$ besides
$A,B,C$ the plane $M$ intersects $\gamma$ at the points $P,Q.$
Therefore, we have the $4$ secant lines $AC,\ BC,\ PC,\ QC$ on $M.$
By repeating this argument for the planes $\langle r\cup D\rangle ,$
$\langle r'\cup A\rangle ,$ $\langle r'\cup B\rangle ,$ we get $16$
spurious secant lines, $4$ of them have been counted twice. Hence,
$I_r\cdot I_{r'}=2.$

It follows
easily that the index of $\{\thinspace I_r\thinspace\}_{\pp{r}}$ is also $2.$
\end{proof}

\medskip

\begin{proof} {\bf of Lemma \ref{famgiusta}}
Let us remark first of all that the curves $\gamma$ and $I_r$ are birational.
Indeed let $r\cap\gamma =\{ A,B\}.$ If $P\in\gamma ,$ and $P\notin r,$
then the plane $\langle r\cup P\rangle$ intersects $\gamma$ at the points
$A,B,P,C,D.$ We get a birational map $f\colon\gamma\to I_r$ by setting
$f\colon P\mapsto  \overline{CD}.$

We fix now a general secant line $r$ of $\gamma.$ Starting
from the just constructed map $f$, we can also construct,
in a canonical way, a  map $f^{(2)}\colon\gamma^{(2)}\to I_r^{(2)},$ which
is again birational.

In the first step of the proof of Lemma \ref{curvabuona} we have constructed a
birational map $\psi\colon\Sigma\cdots\to\sigma^{(2)}.$ Since
 $\gamma$ and $\sigma$ are birational, we get also
a map $\varphi\colon\Sigma\cdots\to\gamma^{(2)}.$

{}Finally, the algebraic system $\{\thinspace I_r\thinspace\}_{\pp{r\in S}}$\
allows us to construct a birational map $\chi\colon I_r^{(2)}\cdots\to S$
as follows. Let $a,b$ be a general pair of secant lines of
$\gamma,$ and assume that each of them  intersects $r.$ By Lemma
\ref{gradoind2} we have $I_a\cdot I_b=2;$ one of these
intersections is $r,$ the other one is, by definition, $\chi(a,b).$

\smallskip

If we compose $\varphi$, \ $f$ and $\chi$ we get the desired map
$\tau\colon\Sigma\cdots\to S.$

\smallskip
It remains to show that $\tau(\sigma(g))=I_{\tau(g)}$.
Consider a curve $\sigma(g)$ such that $g$ intersects $\overline
g.$ It is mapped by $\varphi$ to the curve on $\gamma^{(2)}$ formed by all
the pairs of elements of $\gamma$ containing $g.$ Therefore,
$f^{(2)}\circ\varphi$ sends $\sigma(g)$ to the curve on $I_r^{(2)}$ formed 
by all the pairs of elements of $I_r$ containing $f(g)$,
and clearly $\chi$ maps this last curve to $I_{\tau(g)}.$
\end{proof}

\begin{rmk} \label{uuu}
Note that, if $X$ is one of the threefolds found in this
section with $\mu=3,4$, then the Fano scheme $\Sigma$ 
of $X$ is actually irreducible.
\end{rmk}


\section{Every irreducible component $\Sigma_i$ of 
$\Sigma$\ has $\mu_i=1$}


In this section we assume that the family of lines $\Sigma$ on 
$X$ is reducible and that for every irreducible component $\Sigma_i$ 
of $\Sigma$ we have $\mu_i=1$. 

Note that, from $\mu_i=1$ for all $i$ and from Theorem
\ref{sei}, it follows that $s=\mu\leq 6$.

\bigskip

\noindent
{\bf The case $s =2$.}

\medskip

\begin{prop}\label{s2}
Let $X\subset\sp 4$  be a threefold containing two irreducible families of
lines $\Sigma_i$ ($i=1,2$) both  with $\mu_i =1.$ Assume that
$X$ is not a quadric bundle.
Then $X$ is a threefold of degree $6$ with sectional genus
$\pi=1$, projection of a Fano threefold of $\sp 7$, hyperplane section of
$\sp 2\times\sp 2$ (see
\cite{S}).
\end{prop}

\begin{proof}
If $g_1$ is a fixed line of $\Sigma_1$, then the lines of
$\Sigma_2$ meeting it generate the {\it rational}
ruled surface $\sigma_2(g_1)$
having $g_1$ as simple unisecant. Hence $\Sigma_2$
results to be a rational surface. Similarly for $\Sigma_1$.

There are two possibilities regarding
the algebraic system $\{\sigma_2(g_1)\}_{g_1\in\Sigma_1}$, whose dimension
is  two (because $X$
is not a quadric bundle): either it is already linear, or it can be
embedded in a
larger linear system of curves in $\Sigma_2$, which corresponds to
a linear system of rational ruled surfaces on $X$.
We will prove now that the second case can be excluded.

To this end, we reformulate the problem in a slightly different way. We
consider the rational map
$\phi:X\rightarrow
\sp r:=\sp{H^0(\sigma_2(g_1))^*}$ associated to the  complete linear
system
$\mid~\sigma_2(g_1)\mid$. The map $\phi$ sends a point $p$
to the subsystem formed
by the ruled surfaces passing through $p$. From
$\mu_2=1$, it follows that $\phi$ contracts the lines of $\Sigma_2$, which
are therefore
the fibres of $\phi$. Hence $\phi(X)$ is a surface $S$ 
of degree $d=\sigma_2(g_1)^2$. By an argument similar to that of
Proposition \ref{grado}, we have that
$\deg\sigma_2(g_1)=d+2$.

 The inverse images of the hyperplane sections of $S$ are the
surfaces of
$\mid\sigma_2(g_1)\mid$, so $S$ is a surface with rational hyperplane
sections. We replace now $S$
with a general projection  in $\sp 3$, so we  can apply the theorem of
Kronecker--Castelnuovo and we get only three possibilities:
\begin{enumerate}
\item $S=\sp 2$: in this case the considered algebraic system is already
linear and $d=1$;

\item $S$ is a scroll and  $d>1$;

\item  $S$ is a Steiner surface, projection of a Veronese surface, with
$d=4$.
\end{enumerate}
We have to prove that only the first case happens.  Assume by contradiction that
$S$ is like in
$2.$ or $3.$  Note  that any section of
$S$ with a tangent plane is reducible. If $S$ is a scroll, such a section
is the union of a line $l$ with a plane curve $C$ of degree
$d-1$. Let $\pi$ be the
arithmetic genus of
$C$. The following relation expresses the arithmetic genus of a reducible
plane section of $S$:
$\pi+d-2=0$, so $d=2$, $\pi=0$ and $S$ is a quadric.
Moreover $\deg(\sigma_2(g_1))=4$, so a general ruled surface in the linear
system
$\mid\sigma_2(g_1)\mid$ is a scroll of type $(1,3)$ or $(2,2)$. The  case
$(1,3)$ is excluded
because every surface of the system should  have a unisecant line and our
threefold $X$ contains a family 
of lines of dimension exactly $2$. So  a general scroll of the
system should be of type
$(2,2)$, hence contain a $1$-dimensional  family of conics. 
In this case $X$ contains a
$4$-dimensional family of conics, and  a general hyperplane section $X\cap
H$ of its contains a
$2$-dimensional family of conics. By  the usual argument, $X\cap H$ is a
quadric or a cubic
scroll or a Steiner surface: all three  possibilities are easily excluded.

We assume now  that $S$ is a projection of a Veronese surface. In this case
$\deg\sigma_2(g_1)=6$, so a general ruled surface in the linear system
$\mid\sigma_2(g_1)\mid$ is
a scroll of type $(2,4)$ or $(3,3)$. The reducible plane sections of $S$
are unions of
conics and correspond to reducible ruled surfaces on $X$, unions of two
scrolls of degree three.
Necessarily they are both of type
$(1,2)$ so each of them contains a family of conics of dimension $2$: we
conclude as in the
previous case.

\medskip

So we have proved that for both systems of lines $d=1$, hence
 $\deg \sigma_2(g_1)=\deg \sigma_1(g_2)=3$.
Also the  curves in the Grassmannian $\G 4$ corresponding to these ruled
surfaces have degree $3$.
So the surface $\Sigma_i$ ( for $i=1,2$) contains  a linear system of
dimension two of rational
cubics, with self--intersection one: it defines a birational map from
$\Sigma_i$ to $\sp 2$,
whose inverse map is defined by a linear system of plane cubic curves.
Hence $\deg\Sigma_i\leq
9$ and $\Sigma_i$ has rational or elliptic hyperplane sections.

Moreover there is a natural
birational map between  plane sections of $X$ and some hyperplane sections of
$\Sigma_i$. Precisely, let $H$ be the singular hyperplane section of $\G
4,$ given by lines
meeting a plane $\pi$: then $\Sigma_i\cap H$  represents lines of
$\Sigma_i$  passing through the points of $X\cap\pi$.
Since  there is only one line of $\Sigma_i$ through a general point of $X$,
we get the required
birational map between $\Sigma_i
\cap H$ and $X\cap\pi$.

We conclude that also the plane sections of $X$ are rational or elliptic
curves. In particular a
general hyperplane section of $X$ is a surface of $\sp 3$ with the same
property. The case of
rational sections can be excluded using the Kronecker--Castelnuovo theorem as in
Proposition \ref{mi3}. So a hyperplane section of $X$ is a Del Pezzo
surface and $X$ is a (projection
of) a Fano threefold. Looking at  the list of Fano threefolds we get the
proposition.
\end{proof}

\bigskip

\noindent
{\bf The case $s >2$.}

\medskip

If $\Sigma$ has three or more components, a new situation can appear, precisely
$X$ could be a quadric bundle in more than one way.

{}For example, if $X=\sp 1\times\sp
1\times\sp 1$ (or one of its  projections),  $\Sigma$ has three components
with $\mu_i=1$, so that
there are three lines passing through any point $P$ of $X$, one for each
of the three systems.
The lines of a system
$\Sigma_i$ meeting a fixed line of another system $\Sigma_j$ fill up a
smooth quadric, so the
surfaces
$\sigma_i(g_j)$ are all quadrics. Moreover the $1$-dimensional families
$\{\sigma_i(g_j)\}_{g_j\in\Sigma_j}$ and
$\{\sigma_j(g_i)\}_{g_i\in\Sigma_i}$ coincide. Hence
there are three different structures of quadric bundle on $X$
 giving raise to six families of conics
in  $\G 4$.

\medskip

Let $X$ be a threefold of $\sp 4$ covered by $s\geq 3$  two-dimensional
families of lines
$\Sigma_i$, $i=1,\ldots, s$. We distinguish the following two cases:
\begin{itemize}

\item there exists a pair of indices $(\bar\i,\bar\j)$ such that the
family
$\{\sigma_{\bar\i}(g_{\bar\j})\}_{g_{\bar\j}\in\Sigma_{\bar\j}}$ has
dimension two;

\item for all $(i,j)$, $\dim \{\sigma_i(g_j)\}_{g_j\in\Sigma_j}=1$.

\end{itemize}

In the first case, we consider only the two components
$\Sigma_{\bar\j}$ and $\Sigma_{\bar\i}$: we
can argue on these components as we did in the case $s=2$, obtaining that
$X$ has to be a
projection of a Fano threefold. Since there are no Fano threefolds
satisfying our assumption, we can exclude the first case.

Therefore, if $s\geq 3$, necessarily the surfaces $\sigma_i(g_j)$ are
smooth quadrics for all pair
$(i,j)$. To get the classification, our strategy 
will be the usual one: to fix three of the families
of lines and argue with them. Our result is:

\begin{prop}\label{s3}
Let $X$ be a threefold of $\sp 4$ containing three or more 
irreducible families of lines $\Sigma_i$
all  with $\mu_i =1.$
Then $X$ is a threefold of degree $\leq 6$ with sectional genus
$\pi=1$, projection of $\sp 1\times\sp 1\times\sp 1$.
\end{prop}

\begin{proof}
{}For every pair of indices $(i,j)$ and general $g_j\in\Sigma_j$, the surface
$\sigma_i(g_j)$ is a
smooth quadric and it is clear that the linear systems
$\{\sigma_i(g_j)\}_{g_j\in\Sigma_j}$ and
$\{\sigma_j(g_i)\}_{g_i\in\Sigma_i}$ coincide: we call it $\Sigma_{ij}$. We
want to study  the
intersection of two quadrics belonging to two families of the form
$\Sigma_{ik}$ and $\Sigma_{jk}$,
$i\neq j$.

Let us remark first that, if $g_j$, $g_k$ are two general coplanar lines  in
$\Sigma_j$, $\Sigma_k$ respectively, 
then two cases are possible: either the plane $\langle g_j,
g_k\rangle$ does contain a line of
$\Sigma_i$, or it does not. In the first case $X$ is a cubic (Prop.
\ref{notre}). So if $\deg X>3$ and $p\in\sigma_j(g_k)$,
$p\not\in g_k$, then $p\not\in\sigma_i(g_k)$.
This immediately implies that $\sigma_j(g_k)\cap \sigma_i(g_k)=g_k$.
Let us consider now $\sigma_j(g_k)\cap \sigma_i(g_k')$: it can be written
also as
$\sigma_k(g_j)\cap \sigma_k(g_i^t)$ for a fixed $g_j\in\Sigma_j$ and
$g_i^t$ varying in a ruling of
the second quadric. 
Now $g_j$ certainly meets all the quadrics of
$\Sigma_{ik}$ and is not contained
in any of them, so there exists a $\bar t$ such that $g_j$ and $g_i^{\bar
t}$ meet at a point $q$.
Let $\overline{g_k}$ be the line of $\Sigma_k$ through $q.$ Then:
$$
\sigma_j(g_k)\cap \sigma_i(g_k')=\sigma_k(g_j)\cap 
\sigma_k(g_i^{\bar
t})=\sigma_j(\overline{g_k})\cap
\sigma_i(\overline{g_k}),
$$
so we fall in the previous case. We conclude that two general quadrics of
these families meet along a
line of the family having the common index.

As a consequence, we have that through a general point $p$ of $X$ there 
pass one quadric of the family
$\Sigma_{ij}$ and one line of $\Sigma_k$.

Now, we embed the $\sp 4$ containing $X$ as a subspace of a $\sp 7,$ and
call $Y\subset\sp 7$  the image of the Segre embedding $\sp 1\times\sp 1
\times\sp 1\to\sp 7.$
If $Q\subset X$ is a fixed general quadric of the family $\Sigma_{12}$, by
acting on $Y$
with an element of the projective linear group, we can assume that
$Q\subset Y$ as well.
Let $L\subset\sp 7$ be a linear subspace of dimension $5,$ in \lq\lq general
position" with respect to $X,$ i.e. $L\cap X$ is a curve.
Let $\Sigma_1',\ \Sigma_2'$ and $\Sigma_3'$ denote the three families
of lines on $Y;$ to fix ideas, assume that $Q$ contains lines of the families
 $\Sigma_1',\ \Sigma_2'$ on $Y.$

\smallskip

We define a rational map $\alpha\colon X\backslash 
L\thinspace\cdots\to Y$ as follows. 
Let $p\in X$ be general; then, the line $r\in\Sigma_3$, such that
$p\in r$,
intersects $Q$ at a single point $p'.$ Let $r'\in\Sigma_3'$ be the line (on $Y$)
containing $p'.$ Set
$
\alpha (p)\colon =\langle L\cup p\rangle\cap r'\
$.
It is clear that $\alpha$ is birational. Moreover, by considering the case of
a hyperplane through $L,$ we see that $\alpha$ takes hyperplane sections of
$X$ to hyperplane sections of $Y.$

\smallskip

There are suitable $\sp 3$\thinspace 's in $\sp 7,$ let us call $M$ one of them,
such that the restriction $\beta\colon Y\setminus M\thinspace\cdots\to\sp 3$
of the projection $\sp 7 \setminus M\thinspace\cdots\to\sp 3$ is birational.
The inverse map $\beta^{-1}\colon\sp 3\thinspace\cdots\to Y$ is defined by
a linear system \ $\vert 3H_{\scriptscriptstyle{\sp 3}}-l_1-l_2-l_3\vert ,$
where the $l_i$\thinspace 's are three lines, pairwise skew.

Since $\alpha$ takes hyperplane sections of $X$ to hyperplane sections of $Y,$
the birational map $(\beta\circ\alpha )^{-1}\colon\sp 3
\thinspace\cdots\to X$ is defined by a linear subsystem of \
$\vert 3H_{\scriptscriptstyle{\sp 3}}-l_1-l_2-l_3\vert ,$ i.e. $X$ is a
projection of $Y=\sp 1\times\sp 1\times\sp 1 ,$ and the proof is complete.
\end{proof}

\vspace{1cm}

\begin{flushright}
Dipartimento di Scienze Matematiche\\
Universit\`a di Trieste\\
34127 -- {\sc Trieste}\\
{\it Italia}
\end{flushright}


\bigskip


\begin{thebibliography}{999}


\bibitem{A} {\sc E. Arrondo}: {\it Projections of Grassmannians of lines
and characterization of
Veronese varieties }, J. Algebraic Geom. {\bf 8} (1999), no. 1, 85--101.


\bibitem{B}{\sc M. Baldassarri}:
{\it Le variet\` a pluririgate a tre dimensioni}, Rend. Sem. Mat. Univ. Padova,
{\bf 19} (1950), 172-200.


\bibitem{Bo}{\sc E. Bompiani}: {\it Sulle variet\`a a $k$ dimensioni con
$\infty^k$ rette},
Rend. Acc. Lincei (8) {\bf 1} (1946), 1001-1004

\bibitem{CG}{\sc H. Clemens -- Ph. Griffiths}:
{\it The intermediate Jacobian of the cubic threefold}, Annals of Math.,
{\bf 95}
(1972), 281-356.

\bibitem{C}{\sc F. Conforto}: Le superficie razionali, Zanichelli, Bologna,
1939

\bibitem{Ful}{\sc W. Fulton}:
Intersection theory, Springer -- Verlag, 1984.

\bibitem{GKZ} {\sc I.M. Gelfand -- M.M. Kapranov --A.V. Zelevinsky}: Discriminants,
resultants and multidimensional determinants, Birkh\" auser, 1994

\bibitem{GH2}{\sc Ph. Griffiths -- J. Harris}:
Principles of Algebraic Geometry, John
Wiley \& Sons, 1978.

\bibitem{GH1}{\sc Ph. Griffiths -- J. Harris}:
{\it Algebraic geometry and local differential geometry}, Ann. scient. \' Ec.
Norm. Sup., {\bf 12} (1979), 355-432.

\bibitem{H}{\sc J. Harris}: Algebraic Geometry, Springer -- Verlag 1992.

\bibitem{H2}{\sc J. Harris}: Curves in projective space,
Les Presses de L'Universit\' e de Montr\' eal, 1982.

\bibitem{Kl}{\sc S.L. Kleiman}: {\it The transversality of the general
translate}, Comp. Math., {\bf 28} (1974), 287-297.

\bibitem{L}{\sc J.M. Landsberg}: {\it Is a linear space contained in a
submanifold? On the number of derivatives needed to tell}, 
J. reine angew. Math., {\bf 508} (1999), 53-60.

\bibitem{LP}{\sc A. Lanteri -- M. Palleschi}:
{\it Projective manifolds containing many rational curves}, Indiana Univ. Math.
J., {\bf 36} (1987), 857-865.

\bibitem{N}{\sc A. Nobile}:{\it Variation du genre des
fibres d'une famille de courbes},  Semin. Theor. Nombres, Univ. Bordeaux I
1986/1987, Exp. No.10, 8
p. (1987).

\bibitem{Pe}{\sc C.S. Peterson}: {\it Applications of liaison theory to
schemes supported on lines,
grow of the deficiency module and low rank vector bundles}, Ph. D. Thesis,
Duke University (1994).

\bibitem{Ro}{\sc E. Rogora}:{\it Metodi proiettivi per lo studio di alcune
questioni
relative alle variet\`a immerse}, Ph. D. Thesis, Universit\`a di Roma I (1993)


\bibitem{bS} {\sc B. Segre}: {\it Sulle $V_n$ contenenti pi\`u
di $\infty ^{n-k} S_k$}, {I},  Lincei - Rend. Sc. fis.
mat. e nat. {\bf 5} (1948),  193--197, II,
  Lincei - Rend. Sc. fis.
mat. e nat.  {\bf 5} (1948),  275--280

\bibitem{bS2} {\sc B. Segre}: {\it Bertini forms and Hessian matrices}, J. Lond. Math. Soc.
{\bf 26} (1951), 164--176

\bibitem{S}{\sc M. E. Serpico}: 
{\it Le variet\` a di  Fano sezioni iperpiane di $\sp 2\times\sp
2$}, Le Matematiche,
{\bf 35} (1980), 113-121.

\bibitem{Si}{\sc C. H. Sisam}: 
{\it On Varieties of Three Dimensions with Six Right Lines through
Each Point}, Amer. J. Math.,
{\bf 52} (1930), 607-610.

\bibitem{T}{\sc A. N. Tjurin}: {\it On the intersection of quadrics},
Russian Math. Surveys, {\bf 30} (1975), 51-105.

\bibitem{To} {\sc E. Togliatti}: {\it Sulle variet\`a a $k$ dimensioni
contenenti almeno
$\infty^k$ rette},
Atti Accad. Torino {\bf 57} (1921)



\end{thebibliography}
\end{document}